\documentclass[12pt]{elsart}
\usepackage{amssymb,amsmath,amsfonts,latexsym}

\def\be#1\ee{\begin{equation}#1\end{equation}}
\newcommand{\bea}{\begin{eqnarray}}
\newcommand{\eea}{\end{eqnarray}}
\newcommand{\beaa}{\begin{eqnarray*}}
\newcommand{\eeaa}{\end{eqnarray*}}

\newcommand{\bee}{\begin{enumerate}}
\newcommand{\eee}{\end{enumerate}}

\def\norm#1{\left\|#1\right\|}             

\def\P{{\mathbb{P}}}
\newcommand{\pr}[1]{\P\left\lbrace#1\right\rbrace}
\def\e{{\mathbb{E}}}\def\E{{\mathbb{E}}}
\newcommand{\ind}{1\hspace{-0.098cm}\mathrm{l}}
\let\BFseries\bfseries\def\bfseries{\BFseries\mathversion{bold}}

\def\pr#1{\P\left\lbrace#1\right\rbrace}

\newcommand{\diam}{\operatorname*{diam}}
\def\N{{\mathbb{N}}}
\def\R{{\mathbb{R}}}

\def\tpsi{\widetilde\Psi}
\def\hpsi{\widehat\Psi}

\newcommand{\eps}{\varepsilon}

\def\al{\alpha}

\def\Z{\mathcal Z}
\begin{document}

\begin{frontmatter}
\title{Small Deviation Probability via Chaining}
\author[fa]{Frank\ Aurzada\corauthref{cor}},
\corauth[cor]{Corresponding author.}
\ead{aurzada@math.tu-berlin.de}
\author[ml]{Mikhail\ Lifshits}
\ead{lifts@mail.rcom.ru}
\address[fa]{Technische Universit\"{a}t Berlin, Institut f\"{u}r Mathematik, Sekr.\ 
MA 7-5, Stra\ss e des 17.\ Juni 136, 10623 Berlin, Germany.}
\address[ml]{St. Petersburg State University, 198504 Stary Peterhof, Dept of 
Mathematics and Mechanics, Bibliotechnaya pl., 2, Russia.}

\begin{abstract}
We obtain several extensions of Talagrand's lower bound for the small deviation 
probability using metric entropy. For Gaussian processes, our investigations 
are focused on processes with sub-polynomial and, respectively, exponential 
behaviour of covering numbers. The corresponding results are also proved for 
non-Gaussian symmetric stable processes, both for the cases of critically small 
and critically large entropy. The results extensively use the classical chaining 
technique; at the same time they are meant to explore the limits of this method.
\end{abstract}

\end{frontmatter}

\begin{center} This is the extended version of a paper that is to appear in\\ {\it Stochastic Processes and Their Applications}.\end{center}

\section{Introduction and main results}
\subsection{Motivation}
General small deviation problems attracted much attention
recently due to their deep relations to various mathematical topics
like operator theory, quantization, strong limit laws in
statistics, etc., cf. the surveys\ \cite{LS,Lif}.

The first goal of this article is to extend the well-known Talagrand 
lower bound for the small deviation probability to the case of Gaussian random 
functions with not necessarily regularly varying behaviour of their metric 
entropy.

Before recalling the known results and stating the new ones, let us introduce 
the necessary notation. Consider a centered Gaussian random function 
$X(t),t\in T$, $T\neq \emptyset$, and assume there exists a separable version of $X$ that we 
consider in the sequel. Assume furthermore that the parameter set $T$ equipped 
with quasi-metric $\rho(s,t)^2=\e(X(t)-X(s))^2$, usually referred to as 
Dudley metric, is a relatively compact metric space. Let 
$$
  N(\eps) :=  \min\lbrace n \in \N \,|\, \exists t_1, \ldots, 
  t_n \in T~\forall t\in T~ \exists i~:~ \rho(t,t_i)\leq \eps\rbrace
$$ 
denote the covering numbers of $(T,\rho)$ and $\sigma:=\diam(T)$. Obviously,
$N(\eps)=1$ whenever $\eps\ge\sigma$. Covering numbers present a common 
quantitative measure for the entropy of the space $(T,\rho)$.

At some places we use the following notation for strong and weak asymptotics. 
For two functions $f$ and $g$, $f(x) \sim g(x)$, as $x\to 0$, means that 
$f(x)/g(x) \to 1$, as $x\to 0$. On the other hand, we use the notation 
$f(x) \preceq g(x)$, as $x\to 0$, if  $\limsup_{x\to 0} f(x)/g(x)<\infty$. 
We also write $g(x)\succeq f(x)$ in this case. Furthermore, we write 
$f(x) \approx g(x)$, as $x\to 0$, if $f(x)\preceq g(x)$ and $g(x)\succeq f(x)$. 
The notation is defined analogously for sequences.

Talagrand's lower bound from \cite{Tal}, which became by now classical 
in the form given by M.\ Ledoux \cite[p.\ 257]{Led}, reads as follows.

\begin{thm} 
Assume that $N(\eps)\le \Psi(\eps)$ for all $\eps>0$ and let the 
bound $\Psi$ satisfy the regularity assumptions
 \begin{eqnarray} \label{c121} 
     C_1\, \Psi(\eps)  &\le&  \Psi\left(\frac \eps 2\right), \qquad  \sigma> \eps>0,
 \\  \label{c122}
      \Psi\left(\frac \eps 2\right) &\le& C_2 \, \Psi(\eps),\qquad \eps>0,
 \end{eqnarray}
with some $C_2>C_1>1$.
Then
\be \label{talest}
    \log \pr{ \sup_{s,t\in T} |X(s)-X(t)|\le \eps } \ge - K\Psi(\eps),
    \qquad \eps>0,
\ee 
with $K>0$ depending only on $C_1,C_2$.\label{thm:tal}
\end{thm}

This result works perfectly well and provides sharp estimates 
for many cases where $\Psi$ is a polynomial-type function. Unfortunately, 
on the one hand, it does not apply to slowly varying bounds, e.g. 
$\Psi(\eps)=|\log\eps|^\beta$, since $C_1>1$ in (\ref{c121}) is impossible 
for such functions. Neither is this theorem applicable to exponential bounds, 
e.g.\ $\log \Psi(\eps)=\eps^{-\gamma} |\log\eps|^\beta$, since it is not 
possible to find $C_2<\infty$ in this case.

Moreover, it is easy to see (cf.\ e.g.\ Example~1 below) that in such cases 
the estimate (\ref{talest}) fails in its present form. However, recently, a number of works appeared where small deviations are studied 
for cases with rather arbitrary behaviour of entropy, see e.g. \cite{Lin,LinShi}. 
In particular, a slow increase of $N(\eps)$ when $\eps$ tends to zero is not 
excluded at all. It is therefore desirable to have a version of 
Theorem~\ref{thm:tal} with a wider application range.

The objectives of this article are as follows. Firstly, we show that a more 
careful estimation in the original proof of Talagrand leads to a generally 
applicable lower bound (Theorem~\ref{thm:g}), which, in particular, in the 
case of slow entropy behaviour returns a correct bound.

In the case of large entropy behaviour, we complete the standard approach by 
combining the classical chaining arguments with the use of Laplace transform 
techniques. To the knowledge of the authors, this has not been applied before; 
and it is their belief that the idea could be used successfully in other 
contexts. For this reason, Section~\ref{sec:chaining} is devoted to the 
chaining technique.

Furthermore, the considerations will show that the classical 
chaining idea leads to ``sum of maxima'' type expressions. Namely, classical 
chaining essentially yields estimates of the form 
$$
   \sup_{s,t\in T} |X(s)-X(t)| \leq 
   2 \sum_{k=0}^\infty \eps_k \max_{i=1,\ldots, N_k}|\xi_{k,i}'|,
$$ 
where $\xi_{k,i}'$ are -- not necessarily independent -- standard Gaussian 
random variables, $N_k=N(\eps_{k+1})$, and $(\eps_k)$ is some arbitrary 
decreasing sequence. The above estimate could be called ``uniform'' chaining, 
as opposed to majorizing measure/generic chaining bounds, cf.\ 
\cite{talbook} for a recent description of the theory.

Using the Khatri-\v Sid\'ak inequality allows to replace the $\xi_{k,i}'$ by 
{\it independent} standard Gaussian random variables $\xi_{k,i}$ giving 
\begin{equation} \label{chainingess} 
   \sup_{s,t\in T} |X(s)-X(t)| \leq 2\sum_{k=0}^\infty \eps_k 
   \max_{i=1,\ldots, N_k}|\xi_{k,i}|,
\end{equation} 
where $\leq$ is to be understood in law. The expression on the right-hand side 
is what we will call ``sum of maxima'' type. For the time being, this observation 
has nothing to do with small deviations; note e.g.\ that taking expectations 
of (\ref{chainingess}) immediately yields Dudley's Theorem. However, as we 
demonstrate in this article, a careful estimation of ``sum of maxima'' type terms 
leads to reasonable small deviation results.

Finally, we apply the above-mentioned techniques also to non-Gaussian symmetric 
stable processes, where everything works analogously -- with the natural 
limitations due to the heavy tails. In fact, the most delicate point to be adapted to the non-Gaussian
case is the Khatri-\v Sid\'ak inequality used in the chaining argument. Fortunately, a version
of this inequality for {\it symmetric} stable variables is available, see Lemma~2.1 in \cite{Sam}.

The paper is structured as follows. In Sections~\ref{sec:resgauss} 
and~\ref{sec:resstab} we state the main results of the article, for the cases 
of Gaussian and symmetric $\alpha$-stable random functions, respectively. 
In order to give a taste of the applicability of the results and to present 
the crucial ``sum of maxima'' examples, we consider some important special cases in 
Section~\ref{sec:exa}.

In Section~\ref{sec:chaining}, we recall the classical ``uniform'' chaining 
argument and present the corresponding result for the Laplace transform. 
Section~\ref{sec:smallentropy} contains the proofs of the general estimate, 
which works for slow and polynomial entropy behaviour. The proof is essentially 
the same for Gaussian and symmetric $\alpha$-stable processes. Contrary to this, 
for the large entropy cases, we have to distinguish Gaussian and non-Gaussian 
stable processes, due to their distinct tail behaviour. The proofs in those 
cases are presented in Sections~\ref{sec:largeg} and~\ref{sec:larges}, 
respectively. The article is concluded by some remarks on further extensions 
and related questions in Section~\ref{sec:concl}.

\subsection{The Gaussian case}\label{sec:resgauss}
A version of Talagrand's result that, in particular, includes the case 
of slow increase of entropy is as follows. Let 
\begin{equation} 
   \tpsi(\eps)=
   \int_\eps^{\sigma} \frac{\Psi(u)}{u}\, d u, \qquad 0<\eps\le {\sigma/2},
\label{eqn:fas}\end{equation}
and $\tpsi(\eps)=\Psi(\eps)$ for $\eps\ge \sigma/2$.
We prove the following.

\begin{thm} \label{thm:g} 
Assume that $N(\eps)\le \Psi(\eps)$ for all $\eps>0$ and let the bound $\Psi$
be a non-increasing continuous function satisfying the regularity assumption
\be \label{c2} 
   \Psi\left(\frac \eps 2\right) \le C_2 \, \Psi(\eps),\qquad \eps>0,
\ee
with some $C_2>1$. Then
\be \label{newest}
     \log \pr{ \sup_{s,t\in T} |X(s)-X(t)|\le K_0\eps } \ge - K\tpsi(\eps),
     \qquad \eps>0,
\ee 
with numerical constants $K_0$ and $K>0$, where $K$ depends on $C_2$ and $K_0$ is a universal constant.
\end{thm}

{\bf Comments.}
\medskip

1. We first notice that Theorem~\ref{thm:g} contains Theorem~\ref{thm:tal}. 
Indeed, assumption (\ref{c121}) yields
\[
     \frac {\Psi(u)}{\Psi(\eps)} \le C_1 \left( \frac u\eps\right)^{-h},\qquad 
     \forall \eps\le u,
\]
with $h=\log C_1/\log 2>0 $. We easily obtain from the 
latter inequality that
\be \label{uh}
    \tpsi(\eps)= \int_\eps^{\sigma} \frac{\Psi(u)}{u}\, d u
    \le C_1 \Psi(\eps) \, \eps^{h}\ \int_\eps^\infty u^{-1-h} d u 
    = C_1 h^{-1}\Psi(\eps).
\ee
It is now clear that (\ref{newest}) implies (\ref{talest}).
\medskip

2. Apart from polynomial-type $\Psi$ already covered by Theorem~\ref{thm:tal}, 
the most instructive applications of Theorem~\ref{thm:g} are the following.

a) If $\Psi(\eps)= C |\log\eps|^\beta$ with some $\beta>0$, then 
\[
      \tpsi(\eps)\sim \frac{C}{\beta+1}\, |\log\eps|^{\beta+1}, \qquad 
      \textrm{as}\quad \eps\to 0.
\]
Hence, $N(\eps)\le C |\log\eps|^\beta$ yields
\[ 
   - \log \pr{ \sup_{t,s\in T} |X(t)-X(s)|\le \eps }
   \preceq  |\log\eps|^{\beta+1}, \qquad\textrm{as}\quad \eps\to 0.
\]
\smallskip
b) If $\Psi(\eps)= C \exp\left\lbrace A|\log\eps|^\alpha\right\rbrace$ 
with some $C,A>0$ and $\alpha\in (0,1)$, then 
\[
     \tpsi(\eps)\sim 
     \frac{C}{A}\, |\log\eps|^{1-\alpha} \exp\left\lbrace A|\log\eps|^\alpha\right\rbrace,
     \qquad \textrm{as}\quad \eps\to 0.
\]
Hence, $N(\eps)\le C \exp\left\lbrace A|\log\eps|^\alpha\right\rbrace$ yields
\be \label{explog}
   -\log \pr{ \sup_{t,s\in T} |X(t)-X(s)|\le \eps }
   \preceq |\log\eps|^{1-\alpha} 
   \exp\left\lbrace A|\log\eps|^\alpha\right\rbrace, 
   \ \textrm{as}  \  \eps\to 0.
\ee
\smallskip

We give concrete cases with the above entropy behaviour in Example~1 in 
Section~\ref{sec:exa} below.
\medskip

3. As one can observe from the above-mentioned examples, the ratio of functions 
$\tpsi$ and $\Psi$ ranges between the constant and the logarithmic
function. Actually, this is always true under our assumptions, since for $\eps\le \sigma/2$
\[
\tpsi(\eps)=
\int_\eps^{\sigma} \, \frac{\Psi(u)}{u}\, d u \le  \int_\eps^{\sigma} \, \frac{d u}{u} \, \Psi(\eps)
= \log\frac{\sigma}\eps \ \Psi(\eps)
\]
and 
\be \label{psitpsi}
\tpsi(\eps)=
\int_\eps^{\sigma} \frac{\Psi(u)}{u}\, d u \ge  \int_\eps^{2\eps} \frac{d u}{u}\ \Psi(2\eps)
=\log 2 \ \Psi(2\eps) \ge \frac {\log 2}{C_2}\ \Psi(\eps).
\ee 
\medskip

4. The reader familiar with the theory of Gaussian processes (see e.g.\ \cite{Lif}) 
will surely notice that the integral characteristic $\tpsi$ has much in common 
with the Dudley integral -- the basic entropy tool for the evaluation of 
{\it large} deviations and moduli of continuity of Gaussian processes.
\bigskip

Let us now come to the case of large entropy behaviour. Note that (\ref{c2}) 
restricts  the application range of Theorem~\ref{thm:g} to essentially 
regularly or slowly varying entropy behaviour. However, with the techniques 
presented in this article we can also tackle the case of exponentially 
increasing entropy. One possibility is the following theorem.

\begin{thm} \label{thm:criticalgauss}
Let us assume that 
\begin{equation} \label{eqn:scale}
   \log N(\eps) \leq C \eps^{-\gamma} |\log \eps|^{-\beta}, 
\end{equation} 
with some $0<\gamma<2$ or $\gamma=2$ and $\beta>2$. Then 
$$
    \log \left| \log \pr{\sup_{t,s\in T} |X_t-X_s| \leq \eps}\right| 
    \preceq \eps^{-\frac{2 \gamma}{2-\gamma}} |\log \eps|^{-\frac{2\beta}{2-\gamma}},
    \qquad \text{for $0<\gamma<2$},
$$ 
and 
$$
     \log \log \left| \log \pr{\sup_{t,s\in T} |X_t-X_s| 
     \leq \eps}\right| \preceq \eps^{-\frac{2}{\beta-2}} ,\qquad 
     \text{for $\gamma=2$ and $\beta>2$.}
$$ 
\end{thm}

Note that, due to the classical Dudley Theorem, the above theorem cannot be 
extended beyond $\gamma=2$ and $\beta>2$. Furthermore, it will become 
clear in Examples~3 and~4 that the above bound obtained from (\ref{eqn:scale}) 
cannot be improved by ``uniform'' chaining methods.

\subsection{Stable case}\label{sec:resstab}
Assume now that $X(t),t\in T$, is a symmetric $\al$-stable process, $0<\alpha<2$,
which means that $(X(t_1),\dots, X(t_n))$ is an $n$-dimensional symmetric 
$\al$-stable vector for all choices $t_1,\dots, t_n\in T$, cf.\ \cite{ST}.
We define the quasi-metric related to $X$ by letting $\rho(s,t)$ denote 
the scale parameter of the stable real variable $X(t)-X(s)$; in other words,
$$
    \e \exp\{iu(X(t)-X(s))\} = \exp\{ -\rho(t,s)^\al|u|^\al  \}.
$$
Alternatively, one could choose $(\e|X(t)-X(s)|^r)^{1/r}$ for any fixed 
positive $r<\al$ as a quasi-metric. We assume that, as in the Gaussian case, 
$\sigma:=\diam(T) < \infty$ and $(T,\rho)$ is a relatively compact space.
In what follows, $N(\eps)$ are the covering numbers of the space $(T,\rho)$, 
as defined above.

An analogue of Talagrand's Theorem, i.e.\ our Theorem~\ref{thm:tal}, for the 
stable non-Gaussian case was recently obtained by the first author in 
\cite{Au1}, where it is shown that the result remains true under the additional 
assumption $C_2<2^\al$.
Recall (cf.\ e.g.\ \cite{ST}, p.\ 546) that admitting $C_2>2^\al$ leads to 
processes which may even be not bounded with probability one. Hence there is no 
chance to prove Talagrand's bound for the non-Gaussian case with 
$C_2>2^\al$ in (\ref{c122}). The critical case  $C_2=2^\al$ merits a special 
consideration. It is the case with ``critically large'' entropy, which will 
be handled below.

However, first, we show that Theorem~\ref{thm:g} admits an extension to the stable 
case, too. Namely, the following is true.

\begin{thm} \label{thm:s} 
 Let $X(t),t\in T$, be a symmetric $\al$-stable process, $0<\al<2$. 
 Assume that the corresponding covering numbers satisfy
 $N(\eps)\le \Psi(\eps)$ for all $\eps>0$ and let the bound $\Psi$
 be a non-increasing continuous function satisfying the regularity assumption
 $(\ref{c2})$  with some $1<C_2<2^\al$.
 Then
\[
  \log \pr{ \sup_{s,t\in T} |X(s)-X(t)|\le K_0\eps } \ge - K\tpsi(\eps),
  \qquad \eps>0,
\] 
with a universal constant $K_0>0$, a constant $K>0$ depending only on $\al$ and $C_2$, and where $\tpsi$ is defined in $(\ref{eqn:fas})$.
\end{thm}

The next theorem excludes again slow entropy behaviour but implicitly handles 
the critical case, i.e.\ large entropy behaviour. Let us denote
\[
  \hpsi(\eps)= \int_0^\eps\left(\frac{\Psi(u)}{u}\right)^{\frac{1}{\al+1}} d u.
\]

\begin{thm} \label{thm:sc} 
 Let $X(t)$, $t\in T$, be a symmetric $\al$-stable process, $0<\al<2$. 
 Assume that the corresponding covering numbers satisfy
 $N(\eps)\le \Psi(\eps)$ for all $\eps>0$ and let the bound $\Psi$
 be a non-increasing continuous function satisfying the regularity assumption
\be \label{c1} 
 C_1\, \Psi(\eps)\le \Psi\left(\frac \eps 2\right),\qquad  \sigma \ge \eps>0,
\ee
with some $C_1>1$. Then
\[
\log \pr{ \sup_{s,t\in T} |X(s)-X(t)|\le K_0 \eps } \ge - K\eps^{-\al} \hpsi(\eps)^{\al+1},
\qquad \eps>0,
\] 
with a universal constant $K_0>0$, a constant $K>0$ depending only on $\al$ and $C_1$.
\end{thm}

This theorem also provides a new sufficient condition for the boundedness of stable processes.

\begin{cor}   \label{cor:dudley}
 Let $X(t)$, $t\in T$, be a symmetric $\al$-stable process, $0<\al<2$.  
 Assume that the corresponding covering numbers satisfy 
 $N(\eps)\le \Psi(\eps)$ for all $\eps>0$ and let the bound $\Psi$
 be a non-increasing continuous function satisfying the regularity assumption $(\ref{c1})$. If
\[
  \int_0^\sigma \left(\frac{\Psi(u)}{u}\right)^{\frac{1}{\al+1}} d u < \infty,
\]
 then the process $X$ is a.s.\ bounded.
\end{cor}

Recall that for $0<\al<1$ no sufficient condition for a.s.\ boundedness 
of stable processes in terms of metric entropy had been available so far. 
When $1\leq \al<2$, Theorem~12.2.1 in \cite{ST} provides a sufficient condition, 
which is better than our Corollary~\ref{cor:dudley}, because the integral test 
is slightly weaker and no regularity assumption is required.

We can even go beyond the last theorem in the case $N(\eps)\leq \Psi(\eps)$ with 
$\Psi(\eps):=C \eps^{-\alpha} |\log \eps|^{-\beta}$ with $\beta>0$. 
Note that Theorem~\ref{thm:sc} only works for $\beta>1+\alpha$.

\begin{thm} \label{thm:sc_poly}
Let $N(\eps)\leq C \eps^{-\alpha} |\log \eps|^{-\beta}$ for $\eps<\sigma$. Then
\begin{multline*}
 \log\pr{\sup_{s,t\in T}|X(s)-X(t)|\leq\eps} \\ \geq \begin{cases}
 - K \eps^{-1/(\beta/\alpha-1)} & \max(1,\alpha)<\beta<1+\alpha,\\
 - K \eps^{-\alpha} |\log \eps|^{1+\alpha} & \beta=1+\alpha,\\
 - K \eps^{-\alpha} |\log \eps|^{-\beta+1+\alpha} & \beta>1+\alpha. \end{cases}
\end{multline*} 
\end{thm}

We will show below that these estimates cannot be improved in 
general by the chaining method. In particular, for 
$\beta\leq \max(1,\alpha)$ no estimate can be obtained by uniform chaining. 
It would be interesting to ask what can be done for stable processes using 
majorizing measure/generic chaining techniques.

\begin{rem} \label{rem:dttstable} Similarly to Corollary~$\ref{cor:dudley}$, 
we have that $N(\eps)\leq C \eps^{-\alpha} |\log \eps|^{-\beta}$ with 
$\beta>\max(1,\alpha)$ implies the a.s.\ boundedness of the process. 
Note that this recovers Dudley's Theorem (Theorem~12.2.1 in \cite{ST}) for 
$\alpha\geq 1$ and provides a new Dudley-type theorem for $0<\alpha<1$. 
\end{rem}

\subsection{Some examples} \label{sec:exa}
In the below examples we use, for simplicity, the term symmetric $\alpha$-stable 
for both, the Gaussian ($\alpha=2$) and the non-Gaussian ($0<\alpha<2$) case.

We start with an example that shows that Theorem~\ref{thm:tal} does not 
return the correct bound for slowly varying $\Psi$.

{\bf Example 1 (Logarithmic behaviour of entropy).}
Let $t_n:=2^{-n^{1/\beta}}$ with some $\beta>0$ and let $M$ be an 
independently scattered symmetric $\alpha$-stable random measure on 
$[0,1]$ controlled by the Lebesgue measure. We consider the process 
$$
X_n := M([0,t_n]), \qquad n\geq 1.
$$
It is easy to calculate that $N(\eps)\leq C |\log\eps|^\beta$.

As an example, let us consider $\beta=1$. Note that, if Theorem~\ref{thm:tal} were 
applicable, it would lead to the estimate 
$$
\pr{ \sup_{n,m\geq 1} |X_n-X_m|\le  \eps } \geq C \eps^K,
$$ 
for some $K,C>0$, which is absurd. Instead, we get 
$$
\log \pr{ \sup_{n,m\geq 1} |X_n-X_m|\le K_0 \eps } 
\geq - K |\log \eps|^{2} ,
$$ 
by Theorem~\ref{thm:g} in the Gaussian and Theorem~\ref{thm:s} in the symmetric stable case, which in fact happens to be the correct order.

Analogous arguments give rise to the small deviation behaviour as stated in (\ref{explog}). 
Similar examples (and counterexamples) can be also obtained by using weighted 
sums of independent sequences that are described in Example~2 below.$\Box$

Now we come to the most simple form of symmetric $\alpha$-stable processes, namely, sequences of independent random variables. We investigate what can be said about the small deviations of such sequences in the case of large entropy behaviour.

{\bf Example 2 (Sequence of independent variables).} Let us consider the sto\-chas\-tic process
$X=\left(\sigma_n\xi_n\right)_{n\ge 1}$, where
$\xi_n$ are i.i.d.\ standard symmetric $\al$-stable random variables.

In the Gaussian case, consider the case $\sigma_n \sim (c \log n)^{1/\gamma} (\log \log n)^{-\beta/\gamma}$. Then $\log N(\eps) \leq C \eps^{-\gamma} |\log\eps|^{-\beta}$. Theorem~\ref{thm:criticalgauss} only applies for $\gamma=2$ and $\beta>2$, whereas the problem makes sense even for $\gamma=2$, $\beta=0$, and $c>2$.

In the stable case, the critical situation is obtained when considering $\sigma_n\sim n^{-1/\al}(\log n)^{-\beta/\alpha}$ with $\beta>1$. It is easy to verify that $N(\eps)\approx \eps^{-\al}|\log\eps|^{-\beta}$ and
\begin{equation}
\log\pr{\sup_n |\sigma_n\xi_n |\le \eps}
\approx -\eps^{-\al} |\log\eps|^{1-\beta},\qquad\text{for all $\beta>1$,} \label{eqn:criticalindep}
\end{equation}
cf.\ \cite{Au2}, Section~4.6. Our Theorem~\ref{thm:sc_poly} gives 
weaker results in all cases. In particular, it only works for $\beta>\max(1,\al)$. $\Box$

Let us now come to the crucial ``sum of maxima'' example, that -- as already mentioned in the introduction -- gains its importance as a prototype arising from the chaining estimate.

{\bf Example 3 (Sum of maxima).} Let $\sigma_n>0$ and let 
$N_k\ge 1$ be some integers. Let $(\xi_{k,i}), k,i\ge 1$, be
an array of i.i.d.\ standard symmetric $\al$-stable random variables. Let 
$T=\{ (\ell,s) \in \N^{\infty}\times \{-1,+1\}^{\infty} : 
\ell_k\le N_k, \forall k\ge 1\}$
and set
\[ 
X(\ell,s) = \sum_{k=1}^\infty \sigma_k s_k \xi_{k,\ell_k} \, ,\qquad (\ell,s) \in T.
\]
Note that $X(\ell,s)$ is a symmetric $\alpha$-stable random 
variable with scale parameter
$\left( \sum_k \sigma_k^\alpha \right)^{1/\alpha}$. Then
\begin{equation} \label{ell1_norm2}
S=\sup_{(\ell,s) \in T} |X(\ell,s)| =  \sum_{k=1}^\infty \sigma_k  
\max_{i=1,\ldots, N_k} |\xi_{k,i}|.
\end{equation}
Even if $N_k=1$ for all $k$, we have a nontrivial example of an $\ell_1$-norm,
\begin{equation} \label{ell1_norm}
   \sup_{(\ell,s) \in T} |X(\ell,s)| = \sum_{k=1}^\infty \sigma_k |\xi_{k,1}|. 
\end{equation}

Certain important cases of the ``simplified'' version  
(\ref{ell1_norm}) were studied in \cite{Au2}. We recall 
only one particular case showing that ``simplified''
is not obvious at all. Let $\xi$ be Gaussian and 
$\sigma_k=k^{-1}(\log k)^{-b}$;  then $X$ is bounded for $b>1$ and
\begin{equation} \label{eqn:nochns}
\log\left| \log \pr{ \sup_{t\in T} |X(t)|\le \eps } \right| 
\approx \eps^{-\frac{1}{b-1}},
\end{equation}
while the entropy satisfies 
$\log N(\eps)\approx \eps^{-2}|\log\eps|^{-2b}$ and thus approaches the famous
Dudley-Sudakov border between the bounded and unbounded processes. Our Theorem~\ref{thm:criticalgauss} returns the correct lower bound for (\ref{eqn:nochns}).


As explained in the introduction, this kind of examples provides a sharp power test 
for the chaining method in the small deviation problem.

In the Gaussian case, we obtain the following.
\begin{prop} \label{prop:sumaxgaussian}
Let $S$ be the sum defined in $(\ref{ell1_norm2})$ with $N_k=e^{ 2^{\gamma k} k^{-\beta}}$ and 
$\sigma_k=2^{-k}$ with some $0<\gamma\leq 2$. Then the order given in Theorem~$\ref{thm:criticalgauss}$ is attained for $0<\gamma<2$ or $\gamma=2$ and $\beta>2$, respectively. 
For $\gamma=2$ and $\beta\leq 2$, the process is a.s.\ unbounded. 
\end{prop}

Although formally our theorems cannot be applied here, the considerations in the introduction show that Proposition~\ref{prop:sumaxgaussian} yields the optimality of our theorems in the sense that classical ``uniform'' chaining estimates cannot lead to better estimates.

For the non-Gaussian stable case, we can get the following analog in the respective critical situation.

\begin{prop}\label{prop:summaxstable} 
Let $S$ be the sum defined in $(\ref{ell1_norm2})$
with  $\sigma_k=2^{-k/\alpha} k^{-\beta/\alpha}$ and $N_k=2^k$. Then $S<\infty$
a.s.\ if and only if $\beta> \max(1,\alpha)$ and we have
\[
\log \pr{ S \leq \eps} 
\approx \begin{cases}
 -  \eps^{-1/(\beta/\alpha-1)} & \max(1,\alpha)<\beta<1+\alpha,\\
 -  \eps^{-\alpha} |\log \eps|^{1+\alpha} & \beta=1+\alpha,\\
 -  \eps^{-\alpha} |\log \eps|^{-\beta+1+\alpha} & \beta>1+\alpha. 
     \end{cases}
\]
\end{prop}
$\Box$

Finally, let us consider an example that seems to be closely related to Example~3 and may be important in other circumstances.

{\bf Example 4 (Binary tree).} Let us take an infinite binary tree and 
associate a standard symmetric $\al$-stable 
random variable $\xi_a$ to every edge $a$ of this tree, where we assume all random variables to be 
independent. Let $|a|\ge 1$ denote the level number of an edge $a$. Let $T$ be the set of all finite branches starting from the root of the tree. Furthermore, we take a non-increasing sequence of positive numbers $(\sigma_n)$ 
and consider
\[
X(t)= \sum_{a\in t} \sigma_{|a|}  \xi_a\ , \qquad t\in T.
\]
Then $X(t)$ is a symmetric $\alpha$-stable random variable with scale 
parameter $\left( \sum_{n\le |t|} \sigma_{n}^\alpha \right)^{1/\alpha}$, 
for all $t\in T$, where $|t|$ is the length of the branch.

It is easy to see that this case partially resembles the 
previous example if we set there $N_n=2^n$, although the dependence structures of the two processes are 
substantially different. We have the obvious majoration
\be \label{treemaj}
\sup_{t\in T} |X(t)| \le \sum_{n=1}^\infty \sigma_{n} \, \max_{\{a:|a|=n\}}  |\xi_a|. 
\ee


In the Gaussian case, let us consider the following exemplary situation.

\begin{prop} Let $X$ be the binary tree constructed above with standard normal i.i.d.\ $\xi_a$. \begin{itemize} \item[(a)] Let $\sigma_n=2^{-n/\gamma} n^{-\beta/\gamma}$ with $\gamma>0$ and $\beta\in\R$. Then $$-\log \pr{\sup_{t\in T} \left|X(t)\right| \leq \eps} \approx \eps^{-\gamma} |\log \eps|^{-\beta}.$$
\item[(b)] Let $\sigma_n=n^{-1/2-1/\gamma} (\log n)^{-\beta/\gamma}$ with $0<\gamma<2$ and $\beta\in\R$. Then $$\log\left|\log \pr{\sup_{t\in T} \left|X(t)\right| \leq \eps} \right|\approx \eps^{-\frac{2 \gamma}{2-\gamma}} |\log \eps|^{-\frac{2\beta}{2-\gamma}}.$$ \end{itemize} \label{prop:treesgauss} \end{prop}

The second assertion shows that Theorem~\ref{thm:criticalgauss} cannot be improved since we have $\log N(T,\rho,\eps)\approx \eps^{-\gamma} |\log \eps|^{-\beta}$. However, the method of proof of Proposition~\ref{prop:treesgauss} does not suffice to show the bounds corresponding to the case $\gamma=2$, $\beta>2$ in Theorem~\ref{thm:criticalgauss}. So, there is a gap in the results in this critical case. In fact, it is not clear for which $\beta$ the process is actually bounded when $\gamma=2$; we only know from (\ref{treemaj}) and Proposition~\ref{prop:sumaxgaussian} that $\beta>2$, $\gamma=2$ is sufficient.

For the non-Gaussian stable case, let, in particular, $\sigma_n \sim 2^{-n/\gamma} n^{-\beta/\gamma}$ for some $\gamma>0$, $\beta\in\R$. Then $N(\eps)\approx \eps^{-\gamma} |\log \eps|^{-\beta}$. In this case, we can apply all our theorems. One can also apply the same method used in the proof of Proposition~\ref{prop:treesgauss} to obtain the upper bounds corresponding to Theorem~\ref{thm:tal} for $\gamma>\alpha$:

\begin{prop} Let $X$ be the binary tree constructed above with standard symmetric $\alpha$-stable i.i.d.\ $\xi_a$. Let $\sigma_n=2^{-n/\gamma} n^{-\beta/\gamma}$ with $\gamma>\alpha$ and $\beta\in\R$. Then $$-\log \pr{\sup_{t\in T} \left|X(t)\right| \leq \eps} \approx \eps^{-\gamma} |\log \eps|^{-\beta}.$$ \label{prop:treesstable} \end{prop}

However, the most challenging is the stable non-Gaussian case with $\gamma=\alpha$. In view of (\ref{treemaj}), Proposition~\ref{prop:summaxstable} provides the lower bounds for small deviation probabilities of $X$ whenever $\beta>\max(1,\alpha)$. On the other hand, it is easy to show, by considering the oscillations on each separate level, that $X$ is not bounded when $\beta\le 1$.
Note that, for $\alpha<1$, the process is bounded if and only if $\beta>1$, by Theorem~10.4.2 in \cite{ST}. Observing that
\[ 
   \pr{ \sup_{t\in T} |X(t)|\le \eps/2 } \le
   \pr{ \max_n \sigma_n \sup_{\{a:|a|=n\}} |\xi_a|\le \eps } 
\]
it is easy to show that for any $\beta$
\[
  \log \pr{ \sup_{s,t\in T} |X(s)-X(t)|\le \eps  } 
  \le - K \eps ^{-\alpha}|\log\eps|^{1-\beta}.
\]
There is a gap between this bound and those coming from Proposition~\ref{prop:summaxstable}.
Moreover, we even do not know whether $\beta\in (1,\alpha]$
corresponds to a bounded process $X$. Therefore, many interesting 
questions related to this example remain open.
$\Box$

{\bf Example 5 (L\'evy's Brownian sheet).} Let $\Z$ be a symmetric $\al$-stable 
random measure that is independently scattered on $\R^d_+$ and controlled by the 
Lebesgue measure. For $t\in \R^d_+$ let $[0,t]$ denote the parallelepiped with 
corners $0$ and $t$. Then the random field
\[
Z_\al(t):= \int_{[0,t]} d\Z =\Z([0,t]), \qquad t\in \R^d_+\ ,
\]
is called L\'evy's Brownian sheet. In the Gaussian case this is simply called Brownian sheet. The small deviation problem of $Z_\al$ was studied e.g.\ in \cite{DLKL} for $\alpha=2$ to the end that 
$$\eps^{-2}|\log\eps|^{2d-1} \succeq -\log\pr{\sup_{t\in[0,1]^d} |Z_2(t)|\le \eps}
\succeq \eps^{-2}|\log\eps|^{2d-2},$$ as $\eps\to 0$. For $d=1$, the upper estimate is attained (Brownian motion), whereas, for $d=2$, the lower estimate is the correct one. For $d\geq 3$, the above bounds are the best that are currently known and the true order is unknown. Since $N(\eps)\approx \eps^{-2 d}$, the bound from Theorem~\ref{thm:tal} is far away from being sharp.

In the non-Gaussian case, \cite{LiLin} shows that
\[ -\log\pr{\sup_{t\in[0,1]^d} |Z_\al(t)|\le \eps}\succeq \eps^{-\al} |\log\eps|^{\al(d-1)}, \quad \eps\to 0.\]
For $d>1$, no opposite bound is known. Since $N(\eps)\approx \eps^{-\al d}$, 
neither of our theorems applies to $Z_\al$ for $d\neq 1$. This is just one of many examples where chaining is not an appropriate tool for the evaluation of small deviations.$\Box$

\section{The chaining technique} \label{sec:chaining}
This section is devoted to the basic Dudley-Talagrand chaining argument. For the reader's convenience we shall re-prove it as a separate statement. Following this, we  prove a chaining statement for the corresponding Laplace transform, which turns out to be slightly stronger. However, returning from the Laplace transform to the small deviation probability via Tauberian-type theorems is only possible for regularly varying cases.

These chaining inequalities form the main ingredient of our results. The proofs of our main theorems rely on the following lemmas, appropriate optimization of the parameters in case Lemma~\ref{tallem} is used and appropriate estimates of the involved Laplace transforms if we use Lemma~\ref{talleml}.

\begin{lem} \label{tallem} Let $(\eps_k)_{k\ge 0}$ be a decreasing sequence 
tending to zero such that
$\eps_{0}\ge \sigma$. Let $(b_k)_{k\ge 0}$ be an arbitrary positive sequence. 
Set $b=\sum_{k=0}^\infty b_k$. Then
\be \label{chaining}
\pr{ \sup_{s,t\in T} |X(s)-X(t)|\le 2b} \ge
\prod_{k=0}^\infty \pr{\eps_k |\xi|\le b_k}^{N(\eps_{k+1})}
\ee 
where $\xi$ is a standard normal random variable.
\end{lem}

{\bf Proof.} For any $k\ge 0$, let $T_k$ be a minimal $\eps_k$-net in $T$. Recall that $|T_k|=N(\eps_k)$.
In particular, $|T_{0}|=N(\eps_{0})=1$, since $\eps_{0} \ge \sigma$.

Since $T_{0}$ consists of a single element, we have
\[ \sup_{s,t\in T_{0}} |X(s)-X(t)| = 0,
\]
which provides the induction base. Now we come to the chaining induction step. For any $k\ge 0$,
let $\pi_k:T_{k+1}\to T_k$ be a mapping that satisfies
\[
\max_{t\in T_{k+1}} \rho(t,\pi_k(t))\le \eps_k.
\]
Such a mapping exists by the definition of $T_k$. Then we have the chaining inequality:
for all $s,t\in T_{k+1}$
\begin{multline*}
|X(s)-X(t)| \le |X(s)-X(\pi_k(s))|+|X(\pi_k(s))- X(\pi_k(t))| \\ + |X(\pi_k(t))-X(t)|.
\end{multline*}
Hence,
\[
\sup_{s,t\in T_{k+1}} |X(s)-X(t)|\le 2 \sup_{s\in T_{k+1}}|X(s)-X(\pi_k(s))|+ 
\sup_{s,t\in T_{k}}|X(s)- X(t)|.
\]
By induction, we obtain for any $n\ge 0$,
\begin{equation}\label{eqn:chainingsum}
\sup_{s,t\in T_{n+1}} |X(s)-X(t)|
\le 2 \sum_{k=0}^n \sup_{s\in T_{k+1}}|X(s)-X(\pi_k(s))|.
\end{equation}
Hence, the probability
\[ P_n:=
\pr{
\sup_{s,t\in T_{n+1}} |X(s)-X(t)|\le 2b}
\]
satisfies
\[ P_n
\ge 
\pr{
|X(s)-X(\pi_k(s))| \le b_k,\qquad \forall s\in T_{k+1},  0\le k\le n
}.
\]
By using Khatri-\v Sid\'ak inequality (see e.g.\ \cite[p.\ 260]{Led}) 
and the main property of the mappings $\pi_k$, we get 
\begin{eqnarray*}
\pr{
\sup_{s,t\in T_{n+1}} |X(s)-X(t)|\le 2b}
&\ge& 
\prod_{k=0}^n \prod_{s\in T_{k+1}} \pr{
|X(s)-X(\pi_k(s))| \le b_k }
\cr &\ge& 
\prod_{k=0}^n  \pr{
\eps_k |\xi| \le b_k }^{N(\eps_{k+1})}.
\end{eqnarray*}
Now the assertion follows by a separability argument. $\Box$
\bigskip

Now let us obtain an analog of the chaining lemma, for the corresponding Laplace transform. Recall that it is well-known and has been used at many occasions that considering small deviations of a random variable and the Laplace transform at infinity is equivalent, by the use of Tauberian-type theorems. However, it will turn out that the use of the Laplace transform is technically easier and thus more powerful in a certain sense. In particular, it can be avoided to choose the sequence $(b_k)$, which appears when passing from (\ref{eqn:chainingsum}) to deterministic bounds, which is a somewhat unnecessary step in our context.

\begin{lem} \label{talleml}
Let $(\eps_k)_{k\geq 0}$ be a decreasing sequence tending to zero such 
that $\eps_0\geq \sigma$. Then  
\begin{equation} 
\E \exp\left\lbrace -\lambda \sup_{t,s\in T} |X(t)-X(s)| \right\rbrace \geq  
\prod_{k=0}^\infty \int_0^\infty e^{-y} 
\pr{ 2 \lambda \eps_k |\xi| \leq y}^{N(\eps_{k+1})}\, d y. 
\label{eqn:orig7} 
\end{equation} 
\end{lem}

\noindent{\bf Proof.} By the chaining arguments in the proof of Lemma~\ref{tallem}, we 
obtain (\ref{eqn:chainingsum}). This shows that 
$$
\E e^{-\lambda \sup_{t,s\in T_{n}} |X(t)-X(s)|} \geq 
\E e^{-2 \lambda \sum_{k=0}^n \sup_{s\in T_{k+1}} |X(s)-X(\pi_k(s))|}.
$$
By separability, the left-hand side tends to the Laplace transform we 
wish to evaluate. The right-hand side can be written as 
$$
\int_{\R^{n+1}}  e^{-\sum_{k=0}^n y_k} \pr{ 2 \lambda 
\sup_{s\in T_{k+1}} |X(s)-X(\pi_k(s))| 
\leq y_k, \forall k}\, d (y_0, \ldots, y_{n}).
$$
By the Khatri-\v{S}id\'ak inequality, this is greater or equal to
$$
\int_{\R^{n+1}} e^{- \sum_{k=0}^n y_k} \prod_{k=0}^n 
\prod_{s\in T_{k+1}} \pr{ 2 \lambda |X(s)-X(\pi_k(s))| 
\leq y_k}\, d (y_0, \ldots, y_{n}),
$$ 
which equals
$$
\prod_{k=0}^n \int_{\R} e^{- y_k} \prod_{s\in T_{k+1}} 
\pr{ 2  \lambda |X(s)-X(\pi_k(s))| \leq y_k}\, d y_k.
$$
Note that this is greater or equal to 
$$
\prod_{k=0}^n \int_{\R} e^{- y} \prod_{s\in T_{k+1}} 
\pr{ 2 \lambda \eps_k |\xi| \leq y}\, d y 
= \prod_{k=0}^n \int_{\R} e^{- y} 
\pr{2\lambda \eps_k |\xi| \leq y}^{N(\eps_{k+1})}\, d y,
$$ 
as required in (\ref{eqn:orig7}).$\Box$

\begin{rem} \label{rem:chinggen} 
Let us make an important remark about a slightly 
more general chaining construction. Our calculations still work if we have, 
similarly to $(\ref{eqn:chainingsum})$,
$$
   \sup_{s,t\in T} |X(s)-X(t)|
   \le C\  \sum_{k=0}^\infty \eps_k \max_{i=1,\dots, N_{k+1}} |\xi_{k,i}'|
$$
for some, possibly dependent, standard Gaussian (or, according to the context, 
symmetric stable) variables $\xi_{k,i}'$. In this approach, the $N_n$ are not 
necessarily covering numbers. This observation will be particularly useful 
when considering the tree-based examples. 
\end{rem}

\section{Proofs for the cases with small entropy} \label{sec:smallentropy}
We now assume that the covering numbers admit a reasonable majorant
$\Psi$ and construct, under mildest possible assumptions on $\Psi$, the appropriate
lower bounds for the products appearing in Lemma~\ref{tallem}.

We first show that under (\ref{c2}) the layers with small $\eps_k$ never bring anything really
different from Talagrand's bound.

\begin{lem} \label{smalleps}
Assume that $\Psi$ satisfies $(\ref{c2})$ for $\eps\le\eps_0$.
Then, for any $\eps\in(0,\eps_0)$ and any $r\in(\frac 12, 1)$ it is true that
\be \label{layerl}
\prod_{k=0}^\infty \pr{2^{-k}\eps |\xi|\le r^k\eps }^{N(2^{-k-1}\eps)}
\ge \exp\{-C_3(r)\Psi(\eps)\},
\ee
where $C_3(r)$ depends only on $C_2$ and $r$. 
\end{lem}

{\bf Proof.}\ Since $N(2^{-k-1}\eps)\le \Psi(2^{-k-1}\eps)\le C_2^{k+1}\Psi(\eps)$,
we obviously have
\begin{eqnarray*}
\prod_{k=0}^\infty \pr{2^{-k}\eps |\xi|\le r^k\eps }^{N(2^{-k-1}\eps)}
&\ge&
\prod_{k=0}^\infty \pr{ |\xi|\le (2r)^k }^{C_2^{k+1}\Psi(\eps)}
\cr &=&
\prod_{k=0}^\infty \left( 1- \pr{ |\xi|> (2r)^k }\right)^{C_2^{k+1}\Psi(\eps)}.
\end{eqnarray*}
Since $r>\frac 12$ and $k\ge 0$, we have $\pr{ |\xi|> (2r)^k }\le \pr{ |\xi|> 1 }$,
hence, by using the standard Gaussian tail estimate, we get for some numerical constant $A$,
\begin{equation} 
1- \pr{ |\xi|> (2r)^k } 
\ge \exp\left\lbrace -A   \exp\{ -\frac{1}{2} (2r)^{2k} \} \right\rbrace.  \label{gausstail}
\end{equation}
It follows that
\begin{eqnarray*}
\prod_{k=0}^\infty \pr{2^{-k}\eps |\xi|\le r^k\eps }^{N(2^{-k-1}\eps)}
&\ge& \exp\left\lbrace  -A    \sum_{k=0}^\infty  \exp\{ -\frac{1}{2} (2r)^{2k} \}
C_2^{k+1}\Psi(\eps) \right\rbrace  \\
&=:& 
\exp\left\lbrace  -C_3(r) \Psi(\eps)\right\rbrace,
\end{eqnarray*}
where the sum converges since $r>1/2$. $\Box$
\bigskip

We pass now to the evaluation of the product over the relatively large levels
(small $k$).  Let $r\in (0,1)$ and fix any $\eps>0$. Let $(\eps_k), 0\le k\le n$, 
be a decreasing positive sequence such that $\eps_n=\eps$ and
\be \label{key}
\Psi(\eps_k)\le  r\ \Psi(\eps_{k+1}), \qquad 1\le k\le n-1.
\ee
We set
\[ b_k= r^{n-k}\eps,  \qquad 0\le k<n.
\]

\begin{lem} \label{largeps} With notation introduced above and
 under assumption $(\ref{key})$  we have
\be \label{layer0}
\prod_{k=0}^{n-1} \pr{\eps_k |\xi|\le b_k}^{N(\eps_{k+1})} \ge
\exp\left\lbrace  - C_4(r)\Psi(\eps) -(1-r)^{-1} \, G \ \right\rbrace,
\ee
where $C_4(r)$ depends only on $r$, and
$G =  \sum_{\ell=1}^n   
\left(\log\eps_{\ell-1}-\log\eps_\ell \right) \Psi(\eps_{\ell})
$.
\end{lem}

{\bf Proof.}\ Since for any $k<n$ we have
\[
\frac{b_k}{\eps_k} = r^{n-k}\ \frac{\eps}{\eps_k}= r^{n-k}\ \frac{\eps_n}{\eps_k}\le 1,
\]
it is true that
\be \label{densitybound}
\pr{\eps_k |\xi|\le b_k} \ge c \ \frac{b_k}{\eps_k} 
= c\ r^{n-k}\ \frac{\eps_n}{\eps_k}\ ,
\ee
where $c=(2 \pi e)^{-1/2}$ is a numerical constant. 
On the other hand, it follows from (\ref{key}) that
\be \label{shiftl}
\Psi(\eps_k) \le r^{l-k} \Psi(\eps_\ell), \qquad 1\le k\le \ell,
\ee
in particular,
\be \label{shiftn}
\Psi(\eps_k) \le r^{n-k} \Psi(\eps_n)= r^{n-k} \Psi(\eps), \qquad 1\le k\le n.
\ee
 Therefore,
\[
\prod_{k=0}^{n-1} \pr{\eps_k |\xi|\le b_k}^{N(\eps_{k+1})} \ge
\prod_{k=0}^{n-1} \left(c\ r^{n-k}\ \frac{\eps_n}{\eps_k}\right)^{\Psi(\eps_{k+1})}
= \Pi_1\ \Pi_2,
\]
where
\[
\Pi_1 := \prod_{k=0}^{n-1} \left(c\ r^{n-k} \right)^{\Psi(\eps_{k+1})}
\qquad \textrm{and}\qquad 
\Pi_2 := \prod_{k=0}^{n-1} \left(\frac{\eps_n}{\eps_k}\right)^{\Psi(\eps_{k+1})}.
\]
By using (\ref{shiftn}), we have
\begin{eqnarray*}
|\log\Pi_1| &\le&  \sum_{k=0}^{n-1} \left(|\log c| + |\log r|(n-k) \right)\ \Psi(\eps_{k+1})
\cr  
&\le&  \sum_{k=0}^{n-1} \left(|\log c| + |\log r|(n-k) \right) r^{n-k-1}\ \Psi(\eps)
\cr 
&\le& \sum_{\ell=1}^{\infty} \left(|\log c| + |\log r|\ell \right) r^{\ell-1}\ \Psi(\eps)
=:  C_4(r) \Psi(\eps).
\end{eqnarray*}
Similarly, by using (\ref{shiftl}) and (\ref{shiftn}), we have
\begin{eqnarray*}
|\log\Pi_2| &\le& 
 \sum_{k=0}^{n-1} \left(\log\eps_k-\log\eps_n \right) \ \Psi(\eps_{k+1})
\cr &=&
 \sum_{k=0}^{n-1} \sum_{\ell=k+1}^n \left(\log\eps_{\ell-1}-\log\eps_\ell \right) \ \Psi(\eps_{k+1})
\cr &\leq &
 \sum_{\ell=1}^n   \left(\log\eps_{\ell-1}-\log\eps_\ell \right) \sum_{k=0}^{\ell-1} r^{\ell-k-1}\ \Psi(\eps_{\ell})
\cr &\le&
 (1-r)^{-1} \sum_{\ell=1}^n   \left(\log\eps_{\ell-1}-\log\eps_\ell \right)\ \Psi(\eps_\ell),
\end{eqnarray*}
as claimed above. \ $\Box$
\bigskip

{\bf Proof of Theorem~\ref{thm:g}}.\ 
Let us fix $r\in (1/2,1)$. W.l.o.g.\ $\Psi(\sigma/2)>\Psi(\sigma)$. Therefore, for any $\eps\le \sigma/2$, we can choose $n=n(\eps)\ge 1$ such that 
\[
r^{n-1} \Psi(\eps)> \Psi(\sigma)\ge r^{n}\Psi(\eps).
\]
We choose now the first layer by letting $\eps_0=\sigma$, and the following
$n$ layers from equation
\[ \Psi(\eps_\ell)=r^{n-\ell}\Psi(\eps),\qquad 1\le \ell\le n.
\] 
In particular, we can choose $\eps_n=\eps$.
The choice of $\eps_\ell$\, is possible, since the function $\Psi(\cdot)$ is continuous and 
\[ 
\Psi(\eps)\ge r^{n-\ell}\Psi(\eps)\ge \Psi(\sigma).
\]
Since $\Psi(\cdot)$ is non-increasing,
the sequence $(\eps_\ell)_{0\le\ell\le n}$ is non-increasing as well. 

We put 
\[ b_\ell= r^{n-\ell}\eps,  \qquad 0\le \ell < n,
\]
and apply Lemma~\ref{largeps}. Note that (\ref{key}) is automatically satisfied by the construction of the $\eps_k$. Notice furthermore that for any $1\le\ell\le n$ we have
\[
\Psi(\eps_\ell)\le r^{-1} \Psi(\eps_{\ell-1})
\]
with equality for $2\le \ell \le n$. It follows that
\[
\left( \log\eps_{\ell-1}-\log\eps_\ell \right) \Psi(\eps_{\ell}) \le
r^{-1} \int^{\eps_{\ell-1}}_{\eps_\ell} \frac{d u}u\  \Psi(\eps_{\ell-1})
\le
r^{-1} \int^{\eps_{\ell-1}}_{\eps_\ell} \frac{\Psi(u)}u\, d u .
\]
By summing over $\ell$ we get
\begin{eqnarray*}
 G &=&  \sum_{\ell=1}^n   
 \left(\log\eps_{\ell-1}-\log\eps_\ell \right) \Psi(\eps_{\ell})
 \leq r^{-1} \sum_{\ell=1}^n  \int^{\eps_{\ell-1}}_{\eps_\ell} \frac{\Psi(u)}u \, d u \  
\cr 
 &=& r^{-1}  \int^{\eps_0}_{\eps_n} \frac{\Psi(u)}u \, d u \
  =  r^{-1}  \int_{\eps}^{\sigma} \frac{\Psi(u)}u \, d u,
\end{eqnarray*}
whenever $\eps\le \sigma/2$. 
We obtain from (\ref{layer0})
\be \label{l1}
\prod_{k=0}^{n-1} \pr{\eps_k |\xi|\le b_k}^{N(\eps_{k+1})} \ge
\exp\left\lbrace  - C_4(r) \Psi(\eps) - 
    (1-r)^{-1}r^{-1}   \int_{\eps}^{\sigma} \frac{\Psi(u)}u \, du         
    \right\rbrace.
\ee
We finish the construction by letting
$\eps_{n+k}=2^{-k}\eps$ and $b_{n+k}=r^{k}\eps$ for all positive integers $k$.
By using (\ref{layerl})  we obtain
\be \label{l2}
\prod_{k=0}^\infty \pr{\eps_{n+k} |\xi|\le b_{n+k} }^{N(\eps_{n+k+1})}
\ge \exp\{-C_3(r)\Psi(\eps)\}.
\ee
By plugging (\ref{l1}) and (\ref{l2}) into (\ref{chaining})
and letting $K_0=4\sum_{k=0}^\infty r^k= 4(1-r)^{-1}$
we obtain for
\[
P:=\pr{ \sup_{s,t\in T} |X(s)-X(t)|\le K_0\eps }
\]
that
\be
\label{merge}
P \ge \exp\left\lbrace  - [C_4(r)+ C_3(r) ] \Psi(\eps) - 
    (1-r)^{-1}r^{-1}   \int_{\eps}^{\sigma} \frac{\Psi(u)}u \, d u
    \right\rbrace.
\ee
Finally, consider three cases:
\medskip

a) $0<\eps\le \sigma/2$. Then (\ref{psitpsi}) and (\ref{merge}) yield
\begin{eqnarray*}
P &\ge& \exp\left\lbrace  - [(C_4(r)+ C_3(r))\frac{C_2}{\log 2}+ 
    (1-r)^{-1}r^{-1} ]  \int_{\eps}^{\sigma} \frac{\Psi(u)}u \, d u        
    \right\rbrace \cr
    &=:& \exp\left\lbrace  - C_5 \tpsi(\eps)\right\rbrace.
\end{eqnarray*}

b) $\sigma/2 < \eps\le \sigma$. Then
\[
\int_{\eps}^{\sigma} \frac{\Psi(u)}u \, d u  \le \frac {\Psi(\eps)}{\eps} \cdot \frac \sigma 2 
\le \Psi(\eps),
\]
and hence
\[
P \ge \exp\left\lbrace  - [C_4(r)+ C_3(r)+ 
    (1-r)^{-1}r^{-1} ]\Psi(\eps)        
    \right\rbrace 
    =: \exp\left\lbrace  - C_6 \tpsi(\eps)\right\rbrace.
\]

c) $\eps\ge \sigma/2$. In this case estimate (\ref{layerl}) alone yields
\[
P \ge \exp\left\lbrace  - C_3(r)\Psi(\eps)        
    \right\rbrace 
    = \exp\left\lbrace  - C_3(r) \tpsi(\eps)\right\rbrace.
\]
We choose $K:=\max\{C_3,C_5,C_6\}$ and obtain in all cases
\[
P \ge \exp\left\lbrace  - K \tpsi(\eps)\right\rbrace,
\]
as required.\ 
$\Box$
\bigskip

{\bf Proof of Theorem~\ref{thm:s}}.\ We only indicate here the necessary changes in the proof with respect 
to the Gaussian case. 

The first point is the use of the Khatri-\v Sid\'ak inequality used in the chaining argument. As mentioned in the introduction, this is possible, by Lemma~2.1 in \cite{Sam}. By using this lemma, it was shown in fact in \cite{Au1} (following some ideas of \cite{LifSim}) that the chaining inequality (\ref{chaining}) is still true with the natural replacement of a standard normal random variable $\xi$ by a standard symmetric $\al$-stable random variable.

The second important modification concerns the place where the tail probabilities come into play.
Namely, in Lemma~\ref{smalleps} we must assume that $C_2^{1/\al}/2<r<1$ (recall that $C_2^{1/\al}/2<1$ by 
our theorem's assumption). Instead of (\ref{gausstail}) we have
$$ 1- \pr{ |\xi|> (2r)^k }  \geq \exp\left\lbrace -A\, (2r)^{-\al k} \right\rbrace, $$
where we use the stable tail behaviour:
\begin{equation}
 \pr{ |\xi|\le r  }
\ge \exp\left\lbrace - A r^{-\al} \right\rbrace, \qquad r>0,\label{eqn:stail}
\end{equation}
with some finite positive $A$. Hence this time
\[
\prod_{k=0}^\infty \pr{2^{-k}\eps |\xi|\le r^k\eps }^{N(2^{-k-1}\eps)}
\ge \exp\left\lbrace -C_3 \Psi(\eps)\right\rbrace
\]
where
\[
C_3:= A\,  \sum_{k=0}^\infty (2r)^{-\al k} C_2^{k+1} 
= \frac{A\,C_2}{1- (2r)^{-\al}C_2 }
\]
is finite since $r> C_2^{1/\al}/2$. 

The third point to take care of concerns the density bound used in (\ref{densitybound}). Just note
that the density of a standard non-Gaussian symmetric stable variable is positive and bounded away from zero in any neighborhood of the origin. However, the numerical constant $c$ in (\ref{densitybound}) has to be replaced by the positive number
\[ c:=
\frac{1}{2\pi}\,\int_{-\infty}^{\infty}\cos (u)\, e^{-|u|^\al} d u. 
\]
All other arguments given earlier are valid in the non-Gaussian case, too.\ 
$\Box$
\bigskip

\section{Gaussian case with critically large entropy} \label{sec:largeg}
\subsection{Technical lemmas}
In the following, it will turn out that we have to use a Tauberian-type theorem for 
the Laplace transform that does not seem to be in the literature. The proof is based on, essentially, exponential Chebyshev inequality and a similar estimate. It is in the same spirit as the one for the so-called de Bruijn Tauberian Theorem, i.e.\ Theorem~4.12.9 in~\cite{BGT}, and will therefore be omited.

\begin{lem} \label{lem:taubsuper} 
Let $V$ be a positive random variable. For $\tau>0$ and $\theta\in \R$ 
the following relations are equivalent \begin{eqnarray*}
\log \E e^{-\lambda V} &\,\approx\,& 
- \lambda (\log \lambda)^{-\tau} (\log \log \lambda)^{\theta},
\qquad \lambda\to\infty,
\\ 
\log |\log \pr{ V\leq \eps}| &\,\approx\,&  \eps^{-1/\tau} |\log \eps|^{\theta/\tau}, 
\qquad \eps\to 0. \end{eqnarray*}
Furthermore, let $\theta>0$. Then the following relations are equivalent \begin{eqnarray*}
\log \E e^{-\lambda V} &\,\approx\,& - \lambda (\log \log \lambda)^{-\theta}, 
\qquad \lambda\to\infty,  
\\ 
\log \log |\log \pr{ V\leq \eps}| &\,\approx\,&  \eps^{-1/\theta},
\qquad \eps\to 0.
\end{eqnarray*}
In all statements, the upper (lower) bounds in the assumptions imply lower (upper) bounds in the respective assertions.
\end{lem}

One of the major ingredients of the proofs for the case of critically large entropy is the evaluation of the Laplace transform of the random variable $\max_{i=1,\ldots, N} |\xi_i|$, where $\xi_1, \xi_2, \ldots$ are i.i.d.\ 
standard Gaussian random variables. We start with the case that the argument of the Laplace transform, $L$, is of lower order than $N$.

\begin{lem} \label{lem:gaussiansups1}
Let $\xi_1, \xi_2, \ldots$ be i.i.d.\ standard Gaussian r.v. Then there is a constant $c_1>0$ such that for all $L>0$ and all integers $N\geq 1$ with $2 L\leq N$ we have 
$$  -\log \E e^{-L \max_{i=1,\ldots, N} |\xi_i|} \leq c_1 L \sqrt{\log (N/L)}.$$
Additionally, there is a constant $c_2>0$ such that for all $L\geq 1$ and all 
integers $N\geq 1$ with $2 L\leq N$ we have 
$$  -\log \E e^{-L \max_{i=1,\ldots, N} |\xi_i|} \geq c_2 L \sqrt{\log (N/L)}.$$
\end{lem}

{\bf Proof.} In order to get the first part, note that 
\begin{eqnarray*}
   \E e^{-L \max_{i=1,\ldots, N} |\xi_i|} & =&\int_0^\infty e^{- y} \pr{ L |\xi|\leq y}^N \, d y \\
   &\geq& \int_{L \sqrt{2\log (N/L)}}^\infty e^{- y} \, d y \,\cdot\, 
       \pr{ |\xi|\leq \sqrt{2\log (N/L)}}^N\\
   &=& e^{- L \sqrt{2\log (N/L)}} \, e^{N \log \pr{ |\xi|\leq \sqrt{2\log (N/L)}}}\\
   &\geq& e^{- L \sqrt{2\log (N/L)}} \, e^{- C_1 N \pr{ |\xi|> \sqrt{2\log (N/L)}}}\\
   &\geq& e^{- L \sqrt{2\log (N/L)}} \, e^{- C_2 L} \geq e^{- C_3 L \sqrt{2\log (N/L)}},
\end{eqnarray*}
where we used the assumption $N\geq 2 L$ (steps 5, 6, and 7) and the Gaussian tail (step 6). 

For the reverse inequality note first that 
\begin{multline*} \int_0^\infty e^{- y} \pr{ L |\xi|\leq y}^N \, d y \\ 
\leq \int_0^{L \sqrt{2 \log (N/(2L))} } e^{- y} \pr{ L |\xi|\leq y}^N \, d y 
+ \int_{L \sqrt{2 \log (N/(2L))} }^\infty e^{- y} \, d y. 
\end{multline*}
Here, the second term already admits the required estimate. In order to treat the first term, 
consider the function 
$$
    f(y):= e^{- y} \pr{ L |\xi|\leq y}^N,\qquad y\in \left[0,L \sqrt{2 \log (N/(2L))} \right].
$$
Note that 
$$
   f'(y) = - e^{-y} \pr{ L |\xi|\leq y}^N + e^{-y} N \pr{ L |\xi|\leq y}^{N-1} 
   \phi\left(\frac{y}{L}\right) \, \frac{2}{L},
$$ 
where $\phi$ is the density of the standard normal distribution. Clearly, 
$$
 f'(y) \geq e^{-y} \pr{ L |\xi|\leq y}^{N-1}\left( -1 +  \phi\left(\frac{y}{L}\right) \, 
 \frac{2 N}{L}\right)
$$ 
and 
$$
   \phi\left(\frac{y}{L}\right) \, \frac{2 N}{L} 
   \geq \phi\left(\sqrt{2 \log \frac{N}{2L}}\right) \, \frac{2 N}{L} 
    = \frac{4}{\sqrt{2 \pi}} > 1.
$$ 

Thus, $f$ is increasing and 
\begin{eqnarray*} 
 && \int_0^{L \sqrt{2 \log (2 N/L)} } e^{- y} \pr{ L |\xi|\leq y}^N \, d y \\ 
 &\leq& \int_0^{L \sqrt{2 \log (2 N/L)} }f(L \sqrt{2 \log (2 N/L)})\, d y \\
 &\leq& L \sqrt{2 \log (N/(2L))}\, e^{- L \sqrt{2 \log (N/(2L))} } e^{-c L}\\
 &\leq& e^{\frac{1}{2} \log\left( 2 \log (N/(2L))\right) - L \sqrt{2 \log (N/(2L))} } \\
 &\leq& e^{- \frac{1}{2}  L \sqrt{2 \log (N/(2L))} }, 
\end{eqnarray*}
as long as $L\geq 1$, where we have used that 
\begin{eqnarray*}  
  \pr{ |\xi| \leq \sqrt{2 \log (N/(2L))}}^N &=& 
  e^{N \log \pr{ |\xi|\leq \sqrt{2 \log (N/(2L))}}} 
\\ 
  &\leq&  e^{ - N \pr{ |\xi|>\sqrt{2 \log (N/(2L))}}}   \leq e^{-c L},
\end{eqnarray*}
for some $c>0$. This shows the second assertion. $\Box$
\bigskip

For the sake of completeness, we note that, for very small $L$ we obtain a different behaviour.
 
\begin{lem} \label{lem:gaussiansups1a}
There exist constants $\tilde c_1,\tilde c_2>0$, such that, for all $L\leq 1$ and all integers 
$N\geq 2$, 
$$
   \tilde c_2 L \sqrt{\log N}
   \leq -\log \E e^{-L \max_{i=1,\ldots, N} |\xi_i|} 
   \leq \tilde c_1 L \sqrt{\log N}.
$$  
\end{lem}

{\bf Proof.} Note that 
$$
   \log \E e^{-L \max_{i=1,\ldots, N} |\xi_i|} 
   \approx- L \E \max_{i=1,\ldots, N} |\xi_i| 
   \approx - L\, \sqrt{\log N},
$$ 
by the usual Tauberian-type argument for the Laplace transform at the origin (cf.\ \cite{BGT}) 
and the well-known fact that 
$\E \max_{i=1,\ldots, N} |\xi_i| \approx - \sqrt{\log N}$. 
Here, $\approx$ means that the quotient can be estimated from above and below by 
positive finite constants, which is exactly the assertion.$\Box$ 
\bigskip

The case when $L$ is of larger order than $N$ is as follows.
\begin{lem} \label{lem:gaussiansups2}
Let $\xi_1, \xi_2, \ldots$ be i.i.d.\ standard Gaussian r.v. Then there are constants $c_3, c_4>0$ such that for all integers 
$N\geq 1$ and all $L\geq 2 N$ we have 
$$    c_3 N \log (L/N)
   \leq -\log \E e^{-L \max_{i=1,\ldots, N} |\xi_i|} 
   \leq c_4 N \log (L/N). $$ 
\end{lem}

{\bf Proof.} Note that, for some $c>0$, 
\begin{eqnarray*} 
 && \E e^{-L \max_{i=1,\ldots, N} |\xi_i|} \\
 &=&\int_0^\infty e^{- y} \pr{ L |\xi|\leq y}^N \, d y\\
 &\geq & \int_0^L e^{- y} \pr{ L |\xi|\leq y}^N \, d y
 \geq  \int_0^L e^{- y}  \left(\frac{c y}{L}\right)^N \, d y \\
 &\geq & \left(\frac{c}{L}\right)^N \int_N^{N+1} e^{- y} y^N  \, 
     d y \geq \left(\frac{c N}{L}\right)^N e^{- (N+1)} 
 \geq \left(\frac{c' N}{L}\right)^N.
\end{eqnarray*}
Taking logarithms gives the upper bound. The lower bound is proved 
in the same fashion, namely via using
$$
\int_0^\infty e^{- y} \pr{ L |\xi|\leq y}^N \, d y 
\leq \int_0^L e^{- y} \pr{ L |\xi|\leq y}^N \, d y + \int_L^\infty e^{- y}\, d y.
$$ 
The second term is of lower order, the first term is handled using Stirling's Formula. 
Namely, using the uniform bound for Gaussian density, we see that this term is less than 
$$
   \int_0^L e^{-y} \left(\frac{y}{L}\right)^N\, d y 
   \leq  \left(\frac{1}{L}\right)^N \Gamma(N) 
   \leq  \left(\frac{N}{L}\right)^N,
$$ 
where $\Gamma$ is the Gamma function. Taking logarithms gives the lower bound. $\Box$
\bigskip

The behaviour of the Laplace transfrom is yet different if $L$ is of the same order as $N$.
\begin{lem} \label{lem:gaussiansups3}
Let $\xi_1, \xi_2, \ldots$ be i.i.d.\ standard Gaussian r.v. 
Then there are constants $\tilde c_3, \tilde c_4>0$ such that for all $L>0$ and all $N\in \N$ 
with $L/ 2 \leq N\leq 2L$ we have 
$$     \tilde c_3 L    
    \leq -\log \E e^{-L \max_{i=1,\ldots, N} |\xi_i|} 
    \leq \tilde c_4 L. 
$$  
\end{lem}

The proof is analogous to that of Lemma~\ref{lem:gaussiansups2}.

\subsection{Proof of Theorem~\ref{thm:criticalgauss}}
\noindent {\it Preliminaries:} We use (\ref{eqn:orig7}) with $\eps_k=2^{-k}$. 
This implies that 
$$
   \log \E e^{-\lambda \sup_{t,s\in T} |X(t)-X(s)|} 
   \geq  \sum_{k=0}^\infty \log \E e^{-2 \lambda \eps_k 
   \max_{i=1,\ldots, N(\eps_{k+1})}|\xi_i|}.
$$

Let $\Psi(\eps):=  \exp\left\lbrace C \eps^{-\gamma} |\log \eps|^{-\beta}\right\rbrace$. 
Then, by assumption (\ref{eqn:scale}), 
\begin{equation}   \label{eqn:gcrucial}
   \log \E e^{-\lambda \sup_{t,s\in T} |X(t)-X(s)|} 
   \geq  \sum_{k=0}^\infty \log \E \exp\left\lbrace -2 \lambda \eps_k 
   \max_{i=1,\ldots, \Psi(\eps_{k+1})}|\xi_i|\right\rbrace.
\end{equation}

Let, for the purpose of this proof, $e^r=\lambda$ and 
$$
   F(x):=\log \left(\Psi(2^{-(x+1)}) 2^x\right)
   =2^{\gamma (x+1)} (x+1)^{-\beta} + x \log 2 + \log C.
$$

We split the sum (\ref{eqn:gcrucial}) into three parts: namely, we define
$S_1 :=\sum_{\{ k : \Psi(\eps_{k+1}) \leq  \lambda \eps_k\}}$,
$S_2:=\sum_{\{ k : \Psi(\eps_{k+1})\leq 4 \lambda \eps_k \leq 4 \Psi(\eps_{k+1})\}}$,
and
$S_3:=\sum_{\{ k : \Psi(\eps_{k+1})\geq 4 \lambda \eps_k\}}$.
\medskip

\noindent {\it Evaluation of $S_1$:}  
By Lemma~\ref{lem:gaussiansups2}, it can be estimated from below by 
$$
   -\sum_{\{ k : \Psi(\eps_{k+1}) \leq \lambda \eps_k\}} 
   \Psi(\eps_{k+1}) \log \frac{2 \lambda \eps_k}{\Psi(\eps_{k+1}) } 
   = -\sum_{\{ k : F(k)\leq r \}} \Psi(\eps_{k+1}) (r +\log 2 - F(k)).
$$ 
This can be re-written as 
$$ 
    -\sum_{\{k : F(k) \leq r\}} \sum_{F(k)\leq l \leq r} \Psi(\eps_{k+1})  
    = - \sum_{1\leq l\leq r} \sum_{\,1\leq k \leq F^{-1}(l)}  \Psi(\eps_{k+1}).
$$

It is clear that 
\begin{equation} \label{eqn:notpolynomial} 
   \Psi\left( \frac{x}{2}\right) \geq \Psi(x)^C,
\end{equation} 
for some $C>1$. Using only (\ref{eqn:notpolynomial}) one can show that the inner sum behaves as the 
largest term, which means that the double sum can be estimated from below by
$$
   - c \sum_{1\leq l\leq r}  \Psi(\eps_{F^{-1}(l)+1}) 
   =  - c \sum_{1\leq l\leq r} e^{F(F^{-1}(l)) - F^{-1}(l) \log 2 }.
$$
Using the same argument, this can be estimated again by the largest term in the sum, 
i.e.\ by 
$$
   -c' e^{r - F^{-1}(r) \log 2} .
$$ 
Note that  $F^{-1}(r) \sim \log_2 r^{1/\gamma} + \log_2 (\log r)^{\beta/\gamma}$, 
which shows that the sum $S_1$ behaves, up to a constant, as 
$$
    -\lambda (\log \lambda)^{-1/\gamma} (\log \log \lambda)^{-\beta/\gamma}.
$$

\noindent {\it Evaluation of $S_2$:}   
By Lemma~\ref{lem:gaussiansups3}, it can be estimated by 
\begin{multline*}
   -\sum_{\{ k : \Psi(\eps_{k+1})\leq 4 \lambda \eps_k \leq 4 \Psi(\eps_{k+1})\}} \lambda \eps_k 
   = -\lambda \sum_{\{ k : r \leq F(k) \leq r+\log 4\}} 2^{-k} \\ 
   \geq - \lambda \sum_{\{ k : r \leq F(k) \}} 2^{-k} 
   =  - c \lambda 2^{-F^{-1}(r)}.
\end{multline*} 
This shows that $S_2$ is bounded from below by 
$$
   -\lambda (\log \lambda)^{-1/\gamma} (\log \log \lambda)^{-\beta/\gamma}.
$$

\noindent {\it Evaluation of $S_3$:}  In this case, we can apply the first part 
of Lemma~\ref{lem:gaussiansups1}, which implies that the sum can be estimated by 
$$
   -\sum_{\{ k : \Psi(\eps_{k+1})\geq 4 \lambda \eps_k \}} 
   \lambda \eps_k \sqrt{\log \Psi(\eps_{k+1}) / (2 \lambda \eps_k)}.
$$
Note that this equals 
$$
   -\lambda \sum_{\{ k : F(k)\geq r+\log 4 \}} 2^{-k} \sqrt{F(k)-(r+\log 2)}.
$$
Comparing sum and integral shows that the last term behaves as 
$$
    \approx - \lambda \int_{F^{-1}(r+\log 4)}^\infty 2^{-x} 
    \int_0^{F(x)-(r+\log 2)} y^{-1/2}\, d y\, d x,
$$
which equals
$$   - \lambda \int_0^\infty \int_{F^{-1}(y+r+\log 2)}^\infty 2^{-x}  y^{-1/2}\,  d x\, d y 
     = - 2\lambda \int_0^\infty 2^{-F^{-1}(y+r+\log 2)} y^{-1/2}\, d y.
$$
Recalling that 
$F^{-1}(y) \sim \log_2 y^{1/\gamma} + \log_2 (\log y)^{\beta/\gamma}$ 
shows that the last term behaves as 
$$
   \approx - \lambda \int_0^\infty (y+r+\log 2)^{-1/\gamma} (\log (y+r+\log 2))^{-\beta/\gamma} y^{-1/2}\, d y.
$$ 
Substituting $rz=y$ we obtain 
$$
   \approx -\lambda r^{1/2-1/\gamma} 
   \int_{1}^\infty z^{-1+1/2-1/\gamma} (\log r z)^{-\beta/\gamma} \,d y.
$$ 
Evaluating this, leads to 
$$
   S_3\approx 
   \begin{cases}
   - \lambda (\log\lambda)^{1/2-1/\gamma}  (\log \log\lambda)^{-\beta/\gamma}& 0<\gamma<2\\ 
   -\lambda (\log \log\lambda)^{1-\beta/2}& \gamma=2,\beta>2. 
   \end{cases}
$$ 


Note that the bound for $S_3$ is the dominating term. Applying Lemma~\ref{lem:taubsuper} finishes 
the proof of Theorem~\ref{thm:criticalgauss}.$\Box$

\subsection{Proof of Proposition~\ref{prop:sumaxgaussian}}
The lower bound for the small deviation probability follows, via the observation in Remark~\ref{rem:chinggen} from the proof of Theorem~\ref{thm:criticalgauss}.

For the upper bound, recall that the third sum in the proof of Theorem~\ref{thm:criticalgauss} is the dominating term. If we know that $N\approx\Psi$, all the estimates can be reversed. 
In particular, in order to get an upper bound, we can use the second part of Lemma~\ref{lem:gaussiansups1}, 
by keeping only the sum $\sum_{\{ k : \Psi(\eps_{k+1})\geq 4 \lambda \eps_k \geq 4\}}$.\ $\Box$

\subsection{Proof of Proposition~\ref{prop:treesgauss}}
The lower bound follows from a direct application of Theorems~\ref{thm:tal} and~\ref{thm:criticalgauss}, respectively.

Let us come to the upper bounds. For the sake of readability, we concentrate on (b) and on the special case $\beta=0$, i.e.\ let $\sigma_n=n^{-1/2-1/\gamma}$ for $0<\gamma<2$.

By Anderson's Inequality, cutting the tree into two parts at the root gives: $$\pr{\sup_{t\in T} \left|\sum_{a\in t} \sigma_{|a|}  \xi_{a}\right| \leq \eps} \leq \pr{\sup_{t\in T} \left|\sum_{a\in t, |a|\geq 2} \sigma_{|a|}  \xi_{a}\right| \leq \eps}^2.$$
Iterating the argument yields \begin{equation}\pr{\sup_{t\in T} \left|\sum_{a\in t} \sigma_{|a|}  \xi_{a}\right| \leq \eps}  \leq \pr{\sup_{t\in T} \left|\sum_{a\in t, |a|\geq k+1} \sigma_{|a|}  \xi_{a}\right| \leq \eps}^{2^k}.\label{eqn:c} \end{equation}


We estimate (using a single branch) $$\pr{\sup_{t\in T} \left|\sum_{a\in t, |a|\geq k+1} \sigma_{|a|}  \xi_{a}\right| \leq \eps} \leq \pr{ \left|\sum_{n=k+1}^\infty \sigma_{n}  \xi_{n}\right| \leq \eps},$$ for i.i.d.\ standard normal $(\xi_n)$. This equals in our special case $$\pr{ \left|\left( \sum_{n=k+1}^\infty \sigma_{n}^2\right)^{1/2}  \xi_{0}\right| \leq \eps} \leq \pr{ |\xi_{0}| \leq \frac{\eps (k+1)^{1/\gamma}}{C_\gamma}}.$$

We set $k$ to be the maximal integer such that $k+1\leq K \eps^{-1/(1/\gamma-1/2)}$, with $K$ to be chosen later. Then $$\eps (k+1)^{1/\gamma} \leq \eps \cdot K^{1/\gamma} \eps^{-1/(1-\gamma/2)} = K^{1/\gamma} \eps^{-1/(2/\gamma-1)}\to \infty.$$

Therefore, \begin{multline*}  \log\pr{ |\xi_{0}| \leq \frac{\eps (k+1)^{1/\gamma}}{C_\gamma}}
 = \log \pr{ |\xi_{0}| \leq \frac{K^{1/\gamma} \eps^{-1/(2/\gamma-1)}}{C_\gamma}}
\\
\leq - \pr{ |\xi_{0}| > \frac{K^{1/\gamma} \eps^{-1/(2/\gamma-1)}}{C_\gamma}} 
\leq - \exp\left( - \frac{1}{2} \frac{K^{2/\gamma} \eps^{-1/(1/\gamma-1/2)}}{C_\gamma^2} \right).\end{multline*}

Thus the logarithm of the term in (\ref{eqn:c}) is less or equal to 
$$2^k \, \left(- \exp\left( - \frac{1}{2} \frac{K^{2/\gamma} \eps^{-1/(1/\gamma-1/2)}}{C_\gamma^2} \right)\right) = -\exp\left( k \log 2 - \frac{1}{2} \ldots\right).$$

The term in the exponential equals $$ k \log 2 - \frac{1}{2} \frac{K^{2/\gamma} \eps^{-1/(1/\gamma-1/2)}}{C_\gamma^2} = \left( K \log 2  - \frac{1}{2} \frac{K^{2/\gamma} }{C_\gamma^2}\right) \eps^{-1/(1/\gamma-1/2)} - \log 2.$$ Note that the constant equals $$C':= K \log 2  - \frac{1}{2} \frac{K^{2/\gamma} }{C_\gamma^2} > 0, $$ for $K$ chosen sufficiently small. Thus, $$\log\left( - \log \pr{\sup_{t\in T} \left|\sum_{n=1}^\infty \sigma_{n}  \xi_{t_n}\right| \leq \eps} \right)\geq C' \eps^{-1/(1/\gamma-1/2)} -\log 2,$$ which shows the assertion. The case $\beta\neq 0$ is treated along the same lines (the optimal choice is $k+1\sim K \eps^{-2\gamma/(2-\gamma)} |\log \eps|^{-2\beta/(2-\gamma)}$, with appropriate $K$).

The assertion (a) is proved along the same lines. In fact the proof is even slightly simpler. This time, we have to choose $2^k\sim \eps^{-\gamma}|\log\eps|^{-\beta}$.  $\Box$

\section{Stable case with critically large entropy}\label{sec:larges}
\subsection{Proof of Theorem~\ref{thm:sc}}
Now the construction of small layers from the proof of Theorem~\ref{thm:s} breaks down completely, because
the related evaluation was based on $C_2<2^\al$, which we do not assume anymore. A new construction
is as follows. For $k\ge 0$, let $\eps_k= 2^{-k}\eps$ and 
\[ 
b_k= S^{-1} \left(\eps_k^\al N(\eps_{k+1})\right)^{\frac 1{\al+1}} \eps ,
\]
where
\[ S=S(\eps):= \sum_{k=0}^\infty \left(\eps_k^\al N(\eps_{k+1})\right)^{\frac 1{\al+1}}.
\]
Note that $b=\sum_{k=0}^\infty b_k = \eps$. We use the estimate (\ref{eqn:stail}) which holds for all $r>0$, and obtain
\begin{multline*} \prod_{k=0}^\infty \pr{\eps_k |\xi|\le b_k }^{N(\eps_{k+1})}
\ge \exp\left\lbrace-A \sum_{k=0}^\infty \left(\frac {\eps_k}{b_k} \right)^{\al} N(\eps_{k+1})
\right\rbrace 
\\
= \exp\left\lbrace-A S^{\al} \eps^{-\al} \sum_{k=0}^\infty \eps_k^{\al} 
\left(\eps_k^\al N(\eps_{k+1})\right)^{\frac {-\al}{\al+1}}
N(\eps_{k+1})\right\rbrace = \exp\left\lbrace-A S^{\al+1} \eps^{-\al} \right\rbrace.
\end{multline*}
Now we evaluate $S$. Since $\Psi$ is non-decreasing, we have, for every $k\ge 0$,
\[
\int_{\eps_{k+2}}^{\eps_{k+1}} \left(\frac{\Psi(u)}{u}\right)^{\frac{1}{\al+1}} d u
\ge \Psi\left(\eps_{k+1}\right)^{\frac{1}{\al+1}} 
\int_{\eps_{k+2}}^{\eps_{k+1}} u^{-\frac{1}{\al+1}} d u
= c_\al \left( \Psi\left(\eps_{k+1}\right) \eps_k^\al  \right)^{\frac{1}{\al+1}}.
\]
After summing over $k$, we obtain
\begin{multline*} S \le \sum_{k=0}^\infty \left(\eps_k^\al \Psi(\eps_{k+1})\right)^{\frac 1{\al+1}}
\le c_\al^{-1} \sum_{k=0}^\infty 
\int_{\eps_{k+2}}^{\eps_{k+1}} \left(\frac{\Psi(u)}{u}\right)^{\frac{1}{\al+1}} d u
\\
= c_\al^{-1}
\int_{0}^{\eps_{1}} \left(\frac{\Psi(u)}{u}\right)^{\frac{1}{\al+1}} d u \le
c_\al^{-1} \hpsi(\eps).
\end{multline*}
Therefore,
\[
\prod_{k=0}^\infty \pr{\eps_k |\xi|\le b_k }^{N(\eps_{k+1})}
\ge
\exp\left\lbrace-A c_\al^{-\al-1} \hpsi(\eps)^{\al+1} \eps^{-\al} \right\rbrace.
\]

We do not need to make any changes in the construction and evaluation of higher layers. Therefore,
the estimate (\ref{l1}) remains valid. We just show that both terms from this estimate are dominated
by that of lower layers' bound.

First, we always have for non-increasing $\Psi$,
\[
\hpsi(\eps)^{\al+1} \eps^{-\al} \ge   
\left[ \Psi(\eps)^{\frac{1}{\al+1}}\cdot\eps^{\frac{-1}{\al+1}}\cdot \eps\right]^{\al+1} \eps^{-\al}
=\Psi(\eps).
\]
Second, it follows from (\ref{uh}) that under assumption (\ref{c1})
\[
\tpsi(\eps)= \int_\eps^{\sigma} \frac{\Psi(u)}{u} \, d u
\le C_1 h^{-1}\Psi(\eps),
\]
where $h=\log C_1/\log 2$.

This is enough to get rid of the higher layers. $\Box$
\medskip

{\bf Proof of Corollary~\ref{cor:dudley} and Remark~\ref{rem:dttstable}.} By The\-o\-rem~\ref{thm:sc} and Theorem~\ref{thm:sc_poly}, respectively, it is already clear that the assumptions imply that
\[
\pr{ \sup_{s,t\in T} |X(s)-X(t)|\le K_0 } >0.
\]
Therefore $X$ is bounded with positive probability, which, by the zero-one law in Corollary~9.5.5 in \cite{ST} extends to a.s.\ boundedness. $\Box$
\bigskip

\subsection{Proof of Theorem~\ref{thm:sc_poly}}
We deal with the stable case of critically large entropy, namely when $N(\eps)\leq C \eps^{-\alpha} |\log \eps|^{-\beta}$. The case $\beta>1+\alpha$ is a particular case of Theorem~\ref{thm:sc}. Therefore, let us concentrate on $\max(1,\alpha)<\beta\leq 1+\alpha$.

We are going to use the Laplace technique, i.e.\ Lemma~\ref{talleml} instead of Talagrand's idea from Lemma~\ref{tallem} that was the basis for Theorem~\ref{thm:sc}. Since we deal with a symmetric $\alpha$-stable process we can use the general lower estimate 
(\ref{eqn:stail}). Doing so shows that the term in (\ref{eqn:orig7}) is 
bounded from below by
$$
\prod_{k=0}^\infty \int_0^\infty e^{-y}  \exp\left\lbrace -A  y^{-\alpha} 
(\lambda \eps_k)^\alpha N(\eps_{k+1}) \right\rbrace\, d y.
$$
Using $N(\eps)\leq C \eps^{-\alpha} |\log \eps|^{-\beta}$ and the choice 
$\eps_k=2^{-k}$, we obtain
$$
\prod_{k=1}^\infty \int_0^\infty \exp\left\lbrace -(y +B  y^{-\alpha} 
\lambda^\alpha k^{-\beta}) \right\rbrace\, d y.
$$
We will now need the two following estimates of Laplace integrals, the proofs of which are elemenary and we therefore omit them.

\begin{lem} \label{lem:asympt1}
For $L\to\infty$ we have 
$$\log \int_0^\infty e^{-y-L y^{-\alpha}}\, d y 
\sim - C_\alpha L^{1/(1+\alpha)}.
$$ 
\end{lem}



\begin{lem}  \label{lem:asympt2} 
For $\delta\to 0$ we have 
$$\log \int_0^\infty e^{-y-\delta y^{-\alpha}}\, d y \approx 
\begin{cases} -  \delta^{1/\alpha}& \alpha>1,\\
 - \delta \log 1/\delta &  \alpha=1,\\ 
 - \delta &\alpha<1.\end{cases}
$$ 
\end{lem}



By Lemma~\ref{lem:asympt1} and Lemma~\ref{lem:asympt2} for 
$\beta>\max(1,\alpha)$, $\alpha\neq 1$,
\begin{multline}\sum_{k=1}^{\infty} \log \int_0^\infty 
\exp\left\lbrace -(y +B  y^{-\alpha} \lambda^\alpha k^{-\beta}) \right\rbrace\, d y 
= \sum_{\lambda^\alpha k^{-\beta} >1} + \sum_{\lambda^\alpha k^{-\beta} \leq 1} 
\\ 
\geq -  C_1 \sum_{k<\lambda^{\alpha/\beta}} \lambda^{\alpha/(1+\alpha)} 
k^{-\beta/(1+\alpha)} -  C_2 \sum_{k\geq \lambda^{\alpha/\beta}} 
\lambda^{\min(1,\alpha)} k^{-\beta / \max(1,\alpha)}. 
\label{eqn:lowergenmain}
\end{multline}
For $\max(1,\alpha)<\beta<1+\alpha$, both terms are of order $\lambda^{\alpha/\beta}$. 
This yields that 
$$
\log \E e^{-\lambda \sup_{t,s\in T} |X(t)-X(s)|} 
\geq -C \lambda^{\alpha/\beta}.
$$
By the usual Tauberian-type argument (the so-called de Bruijn Tauberian Theorem, i.e.\ Theorem~4.12.9 in~\cite{BGT}), this shows the assertion 
for the range $\max(1,\alpha)<\beta<\alpha+1$. The argument for $\alpha=1$ is similar.

For $\beta=\alpha+1$, the first term in (\ref{eqn:lowergenmain}) contains an additional logarithm, whereas the second does not and is thus of lower order. This yields $$\log \E e^{-\lambda \sup_{t,s\in T} |X(t)-X(s)|} \geq -C \lambda^{\alpha/(1+\alpha)} \log \lambda,$$ and once again the standard Tau\-ber\-ian-type argument proves the theorem's assertion.$\Box$

\subsection{Proof of Proposition~\ref{prop:summaxstable}}
Recall that we consider the sum of maxima example (Example~3) with 
$\sigma_n = 2^{-n/\alpha} n^{-\beta/\alpha}$ and $N_n=2^n$. 

The lower bound for the small deviation probability follows, via the observation in Remark~\ref{rem:chinggen} applied to $N_n=2^n, \eps_n = 2^{-n/\alpha} n^{-\beta/\alpha}$, from the proof of Theorem~\ref{thm:sc_poly}.
\medskip

\noindent{\bf Proof of the upper bound.} Consider the corresponding Laplace transform
\begin{multline}
\E e^{-\lambda \sum_{n=1}^N \sigma_n 
\max_{k=1,\ldots, N_n} |\xi_{n,k}|}
= \prod_{n=1}^N \int_{\R} e^{- y} 
\pr{ \lambda \sigma_n \max_{k=1,\ldots, N_n} |\xi_{n,k}| 
\leq y}\, d y
\\ = \prod_{n=1}^N \int_{0}^\infty  e^{-y} 
\pr{ \lambda \sigma_n |\xi| \leq y}^{N_n}\, d y. \label{eqn:sumsmaxlp}
\end{multline}

We estimate this term using that $t:=\pr{ |\xi| \leq 1} < 1$ and the equivalent 
to (\ref{eqn:stail}) for large arguments as follows
\begin{multline} 
\int_{0}^\infty  e^{-y} \pr{ |\xi| \leq \frac{y}{\lambda \sigma_n}
}^{N_n}\, d y 
= \int_{0}^{\lambda \sigma_n} + \int_{\lambda \sigma_n}^\infty
\\ 
\leq \int_{0}^{\lambda \sigma_n}  e^{-y} t^{N_n}\, d y + 
\int_{\lambda \sigma_n}^\infty  e^{-y - A \lambda^\alpha 
n^{-\beta} y^{-\alpha}} \, d y. \label{eqn:tailanddens}
\end{multline}

\paragraph*{The case $\beta\geq 1+\alpha$.} Let 
\begin{eqnarray*} 
n \in A_\lambda &:=& \left\lbrace k ~:~
\lambda \sigma_k \leq \lambda^{\alpha/(1+\alpha)} k^{-\beta/(1+\alpha)}, 
k\leq \lambda^{\alpha/\beta}  \right\rbrace 
\\ 
&\supseteq& \left\lbrace k ~:~ \frac{\alpha}{1+\alpha}\, 
\log \lambda \leq k \leq \lambda^{\alpha/\beta} \right\rbrace. 
\end{eqnarray*}
Then the first term in the sum in (\ref{eqn:tailanddens}) can be estimated by 
$$ 
t^{N_n}\leq e^{-C 2^n} \leq e^{-C' \lambda^{\alpha/(1+\alpha)} 
n^{-\beta/(1+\alpha)}}. 
$$ 
On the other hand, the second term in (\ref{eqn:tailanddens}) is less than 
$$
\int_0^\infty e^{-y - A \lambda^{\alpha} n^{-\beta} y^{-\alpha}}\, d y 
\leq e^{- C \lambda^{\alpha/(1+\alpha)} n^{-\beta/(1+\alpha)}},
$$ 
by Lemma~\ref{lem:asympt1} and the fact that $n\leq \lambda^{\alpha/\beta}$. 
Using these estimates, (\ref{eqn:sumsmaxlp}), and letting $N$ tend to infinity, 
we obtain

\begin{multline*} 
\log \E e^{-\lambda S
} 
\leq \sum_{n\in A_\lambda} - C \lambda^{\alpha/(1+\alpha)} n^{-\beta/(1+\alpha)} 
\\ 
\leq - C'  \lambda^{\alpha/(1+\alpha)} \sum_{\frac{\alpha}{1+\alpha}\,  
\log \lambda\leq n \leq \lambda^{\alpha/\beta}} n^{-\beta/(1+\alpha)}. 
\end{multline*}
Note that this term is less or equal to 
$$
\begin{cases} - C''  \lambda^{\alpha/(1+\alpha)} (\log \lambda)^{1-\beta/(1+\alpha)} 
& \beta>1+\alpha,\\ - C''  \lambda^{\alpha/(1+\alpha)} 
(\log \lambda) & \beta=1+\alpha,\end{cases}
$$
which, by the de Bruijn Tauberian Theorem (cf.\ Theorem~4.12.9 in~\cite{BGT}), implies the assertion.

\paragraph*{The case $\max(1,\alpha)<\beta<1+\alpha$.} We let $n \in B_\lambda$, where $B_\lambda := \left\lbrace k ~:~ \lambda^{\alpha/\beta} \leq k \leq 2 \lambda^{\alpha/\beta} \right\rbrace$. Then the second term in (\ref{eqn:tailanddens}) is bounded by 
$$\int_0^\infty e^{-y - A \lambda^\alpha n^{-\beta} y^{-\alpha}} \, d y.
$$ 
Since $A \lambda^\alpha n^{-\beta} \leq A$, we have by Lemma~\ref{lem:asympt2}, 
that the last term is bounded by 
\begin{equation}
 \begin{cases} e^{- C \lambda n^{-\beta/\alpha}}& \alpha>1,
 \\ 
 e^{- C \lambda n^{-\beta} \log( \lambda^{-\alpha} n^{\beta})} 
 &  \alpha=1,\\ e^{- C \lambda^\alpha n^{-\beta}} 
 & \alpha<1.\end{cases} 
 \label{eqn:mtscc}
\end{equation}

On the other hand, the first term in (\ref{eqn:tailanddens}) is bounded by 
$e^{- C 2^n}$, which is certainly smaller than (\ref{eqn:mtscc}). 
Using this, (\ref{eqn:mtscc}), and (\ref{eqn:sumsmaxlp}) and letting $N$ 
tend to infinity we obtain
$$
\log \E e^{-\lambda S
\sum_{n=1}^\infty \sigma_n \max_{k=1,\ldots, N_n} |\xi_{n,k}|
} 
\leq  - C \lambda^{\alpha/\beta},
$$ 
in all three cases. By the de Bruijn Tauberian Theorem (cf.\ Theorem~4.12.9 in~\cite{BGT}), this implies the assertion. 

\paragraph*{The case $\beta\leq \max(1,\alpha)$.} Here we use Kolmogorov's Three 
Series Theorem to show that $S$ is infinite a.s.\ On the one hand, 
it is necessary for the convergence of $S$
that 
$$
\sum_n \pr{ \sigma_n \max_{k=1,\ldots, N_n} |\xi_{n,k}| > 1} < \infty.
$$ 
Using the tail estimate (\ref{eqn:stail}), it is easy to see that this is true 
if and only if $\sum_n \sigma_n^{\alpha} 2^n < \infty$, 
which is violated for $\beta\leq 1$. Thus we are finished for $0<\alpha\leq 1$.

On the other hand, it is necessary for $S$ to be a.s.\ finite that 
\be \label{cut_exp}
\sum_n \E \sigma_n \max_{k=1,\ldots, N_n} |\xi_{n,k}| 
\, \ind_{\{ \sigma_n \max_{k=1,\ldots, N_n} |\xi_{n,k}| \leq 1\}} < \infty.
\ee
Let $\alpha>1$. Note that 
$$\sigma_n\max_{k} |\xi_{n,k}|\, \ind_{\{ \sigma_n \max_{k} |\xi_{n,k}| \leq 1\}} 
=  \sigma_n \max_{k} |\xi_{n,k}| - \sigma_n \max_{k} |\xi_{n,k}| \, 
\ind_{\{ \sigma_n \max_{k} |\xi_{n,k}| > 1\}}.$$
It is easy to show using the tail estimate (\ref{eqn:stail}) that
$$\sum_n \E \sigma_n \max_{k=1,\ldots, N_n} |\xi_{n,k}| 
\ind_{\{ \sigma_n \max_{k=1,\ldots, N_n} |\xi_{n,k}| > 1\}} < \infty 
\quad
\Leftrightarrow\quad \beta>1\ \textrm{and}\ \alpha>1.
$$

However, for $\alpha>1$, $\E \max_{k=1,\ldots, N_n} |\xi_{n,k}| \approx N_n^{1/\alpha}$ 
(cf.\ e.g.\ \cite{MP}, p.\ 271), 
which shows that 
$$
\sum_n \E \sigma_n \max_{k} |\xi_{n,k}| < \infty 
\qquad \Leftrightarrow\qquad \beta>\alpha.
$$ 
It follows that the series (\ref{cut_exp}) diverges when $1<\beta\le\alpha$.
This finishes the proof of Proposition~\ref{prop:summaxstable}. 
$\Box$

\section{Concluding remarks} \label{sec:concl}
{\bf 1.} There is another type of processes with slowly vanishing small 
deviation probabilities. Take for example a stationary Gaussian process 
$X(t), t\in \R$, with quickly decreasing spectral density $f$, say
\[ 
   f(\lambda) = \exp\{-\lambda^2\},\qquad \lambda\in \R.
\]
Then the small deviation probability is vanishing too slowly, e.g.
\[ \lim_{\eps\to 0} \eps^{h}\ \log \pr{ \sup_{s,t\in [0,1]} |X(s)-X(t)|\le \eps } =0,
\qquad \forall \, h>0,
\]
while the covering numbers grow polynomially. Namely, $N(\eps)\approx \eps^{-1}$, 
due to the smoothness of $X$. Such kind of small deviation behaviour can not 
be obtained from our results. It is rather related with extremely good 
approximation of the analytical process $X$ by finite rank processes. 
See \cite{VZ}, for more details and statistical applications.
\bigskip

{\bf 2.} There exists a surprising relation between the small deviations in the 
critical stable and critical Gaussian case, as the following example shows. 
Let $(\xi_n)$ be i.i.d.\ standard Gaussian random variables and let 
$A_n$ be i.i.d.\ totally skewed positive $\alpha/2$-stable random variables. 
Then $\theta_n=A_n^{1/2} \xi_n$ are i.i.d.\ symmetric $\alpha$-stable random variables. 
Let $(\sigma_n)$ be a positive sequence of real numbers that is regularly 
varying for $n\to\infty$ with negative exponent. Then the studies of small 
deviation probabilities $\pr{ \sum_n |\sigma_n \xi_n|^\alpha \leq \eps^\alpha}$ 
and $\pr{\sum_n |\sigma_n \theta_n|^2 \leq \eps^2}$ can be completely reduced 
to each other (at least, on the logarithmic level), by using the Laplace 
transform technique.

In particular, the critical stable case, with 
$\sigma_n\sim n^{-1/\alpha} (\log n)^{-\beta/\alpha}$ considered in 
(\ref{eqn:criticalindep}) with entropy 
$N(\eps) \approx \eps^{-\alpha} |\log \eps|^{-\beta}$, corresponds 
to the Gaussian case with large entropy 
$\log N(\eps) \approx \eps^{-2} |\log \eps|^{-2 \beta/\alpha}$.

Both, the stable and the Gaussian process, are bounded if and only if $\beta>1$.
\bigskip

{\bf 3.} {$\R^d$-valued Processes.}
Let us consider $(X(t))_{t\in T}$ to be a Gaussian or symmetric $\alpha$-stable 
process with values in $\R^d$. Then we define for a Gaussian process 
the analogue to the Dudley metric by 
$$\rho(t,s) := \left(\E \norm{X(t)-X(s)}^2\right)^{1/2},$$ 
replaced by the $r$-th moment for the stable case. 
Here, $\norm{.}$ denotes any norm on $\R^d$. As above we consider the covering numbers 
$N(\eps)$ of the quasi-metric space $(T,\rho)$, which we assume to be relatively compact. 

\begin{prop} 
All the above theorems and corollaries hold literally for the case of an 
$\R^d$-valued Gaussian or symmetric $\alpha$-stable process, respectively. 
\end{prop}
\bigskip

{\bf 4.} {Supremum vs.\ supremum of increments.}
We have formulated all our estimates for the small deviation probability of $\sup_{t,s\in T} |X(t)-X(s)|$. Regarding our results there is no difference to the small ball problem for $\sup_{t\in T} |X(t)|$. This can be seen simply by adding a point $t_0\notin T$ into a new set $T':=T\cup \{t_0\}$ and setting $X(t_0)=0$. Then $$\sup_{t,s\in T} |X(t)-X(s)| \leq 2 \sup_{t\in T} |X(t)| = 2\sup_{t\in T} |X(t)-X(t_0)| \leq 2\sup_{t,s\in T'} |X(t)-X(s)|$$ and $$N(T',\rho,\eps)\geq N(T,\rho,\eps)\geq N(T',\rho,\eps)+1.$$
\bigskip

{\bf Acknowledgements.} The authors are very grateful to W.\ Linde and Z.\ Shi for 
motivating discussions. The research of the first-mentioned author was supported by the DFG Research Center \textsc{Matheon} ``Mathematics for key technologies'' in Berlin. The work of the second author 
was supported by the Edinburgh International 
Institute for Mathematical Sciences (conference ``Metric Entropy and 
Applications in Analysis, Learning Theory and Probability'', 2006) as well as by 
the grants RFBR-DFG 04-01-04000 and NSh. 4222.2006.1.

{
}
\end{document}